\documentclass{amsart}
\usepackage{amssymb,amsmath,amsthm}
\usepackage{tikz}
\usepackage{hyperref}
\hypersetup{
  pdftitle={OpenAI's proof of the local smoothing conjecture: An exposition for harmonic analysts},
  colorlinks=false,
  pdfborder={0 0 1},
  linkbordercolor={1 0 0},
  citebordercolor={0 0.65 0},
  urlbordercolor={0 0.65 0}
}

\numberwithin{equation}{section}
\newtheorem{theorem}{Theorem}[section]
\newtheorem{lemma}[theorem]{Lemma}
\newtheorem{proposition}[theorem]{Proposition}
\theoremstyle{remark}
\newtheorem*{remark}{Remark}

\newcommand{\step}[1]{\par\smallskip\noindent\emph{#1}}

\begin{document}

\title[OpenAI Local Smoothing]{OpenAI's proof of the local smoothing conjecture: A harmonic analytic perspective}

\author{Xiangyu Wang}
\address{Department of Mathematics, University of Illinois Urbana-Champaign}
\email{xw70@illinois.edu}

\begin{abstract}
This note presents OpenAI's proof of the Euclidean local smoothing conjecture for the wave equation. It is written for harmonic analysts who may be new to quaternionic Monge--Ampère equations. We introduce the quaternionic Hessian and explain how the associated determinant equation leads to the critical annular estimate.
\end{abstract}
\maketitle

\section{Introduction}

A solution of the wave equation can be much larger in \(L^p\), \(p>2\), at a later time than it was at time zero: a wave that starts out spread over a neighborhood of a sphere can focus near the center. Fixed-time \(L^p\) estimates therefore lose derivatives. The local smoothing problem asks how much of this loss is recovered when the \(L^p\) norm is averaged over a time interval. Sogge~\cite[Section~4]{sogge} formulated the Euclidean conjecture in his work on circular maximal functions. The connection is that, at high frequency, the Fourier transform of arclength measure on a circle of radius \(t\) has the two leading phases \(\pm t|\xi|\), which are the phases of the half-wave propagators.

For a Schwartz function \(f\) on \(\mathbb R^n\) we use the convention
\(\widehat f(\xi)=\int_{\mathbb R^n}e^{-ix\cdot\xi}f(x)\,dx\),
so that the two half-wave propagators are
\[
    e^{\pm it\sqrt{-\Delta}}f(x)
      :=\frac{1}{(2\pi)^n}\int_{\mathbb R^n}
        e^{i(x\cdot\xi\pm t|\xi|)}\widehat f(\xi)\,d\xi.
\]
We work with the minus sign throughout, so the phase is \(x\cdot\xi-t|\xi|\).
The same estimates hold for the plus sign because
\(e^{it\sqrt{-\Delta}}f=\overline{e^{-it\sqrt{-\Delta}}\overline f}\).

For \(2\le p<\infty\), Peral's fixed-time estimate~\cite{peral} requires \((n-1)(1/2-1/p)\) derivatives. Sogge's conjecture predicts that averaging in time gains almost \(1/p\) derivatives when \(p\ge p_c=2n/(n-1)\). To state the result, for \(\varepsilon>0\) set
\[
  \sigma_p:=\max\left\{0,\frac{n-1}{2}-\frac{n}{p}\right\},
  \qquad s_p:=\sigma_p+\varepsilon ,
\]
and use the inhomogeneous Sobolev norm \(\|f\|_{W^{s,p}}=\|(1-\Delta)^{s/2}f\|_{L^p}\).

\begin{theorem}[Euclidean local smoothing]\label{thm:local-smoothing}
Let \(n\ge 2\), \(2\le p<\infty\), and \(\varepsilon>0\). For every Schwartz function \(f\) on \(\mathbb R^n\),
\[
  \big\|e^{-it\sqrt{-\Delta}}f\big\|_{L^p(\mathbb R^n\times[1,2])}
  \lesssim_{n,p,\varepsilon}
  \|f\|_{W^{s_p,p}(\mathbb R^n)}.
\]
At \(p=2\), one may take \(\varepsilon=0\).
\end{theorem}

For \(2<p\le p_c\) the theorem asks for only \(\varepsilon\) derivatives. For \(p>p_c\) we have \(\sigma_p=(n-1)(1/2-1/p)-1/p\), so averaging in time gains \(1/p-\varepsilon\) derivatives over the fixed-time estimate.

The Euclidean local smoothing conjecture was long regarded as one of the most important and difficult problems in harmonic analysis. It implies the Bochner--Riesz conjecture~\cite[Proposition~39]{beltran-hickman-sogge}, as Sogge already observed in the plane~\cite[(4.5)--(4.6)]{sogge}. Tao~\cite{tao-bochner-riesz} proved that the Bochner--Riesz conjecture implies the Fourier restriction conjecture, which in turn implies the Kakeya conjecture by the work of Beckner, Carbery, Semmes and Soria~\cite{beckner-carbery-semmes-soria}.

\subsection{Earlier work on local smoothing}\label{sec:earlier}

Wolff~\cite{wolff} split a unit-thick neighborhood of the cone at frequency \(R\) into angular sectors of width \(R^{-1/2}\) and used induction on scales to bound the \(L^p\) norm of their sum by \(R^{1/2-2/p+\varepsilon}\) times the \(\ell^p\) sum of the \(L^p\) norms of the pieces. This gave every gain below \(1/p\) for \(p>74\) in two spatial dimensions. Łaba and Wolff~\cite{laba-wolff} extended the method to higher dimensions. The sharp \(\ell^2\) cone decoupling theorem of Bourgain and Demeter~\cite{bdec} then gave the conjectured gain for \(p\ge 2(n+1)/(n-1)\) in every dimension.

In two spatial dimensions, Guth, Wang and Zhang~\cite{gwz} reached the Euclidean critical exponent \(p_c=4\) by proving the reverse \(L^4\) square-function estimate for the cone in \(\mathbb R^3\), which controls a sum of waves moving in different directions by their spacetime square function. Combined with the bound for the square function in terms of the initial data, it resolves the two-dimensional conjecture.

Gan, He, Li and Wu~\cite{ghlw} reached the sharp local smoothing range for compact Riemannian manifolds of odd dimension \(n\ge3\). They combined a wave packet density adapted to Lorentz rescaling with refined decoupling, broad estimates and a two-ends reduction. In even dimensions, the author~\cite{wang-density} refined their approach using a transverse covering \(L^2\) estimate. Together these results give every gain below \(1/p\) for Fourier integral operators with positive conic curvature, including the half-wave propagator on compact Riemannian manifolds, in the range
\[p\ge
  \begin{cases}
    2+\dfrac{8}{3n-3}, & n\ge3\text{ odd},\\[6pt]
    2+\dfrac{8}{3n-2-4/(n+1)}, & n\ge4\text{ even}.
  \end{cases}
\]

\subsection{The critical annular estimate}\label{sec:annular-idea}

For \(\phi\in C_c^\infty(\mathbb R^n\setminus\{0\})\) and \(R\ge 2\), let \(\phi_R\) be the spatial Fourier multiplier defined by
\[
  \widehat{\phi_R f}(\xi)=\phi(\xi/R)\widehat f(\xi).
\]

\begin{theorem}[The critical annular estimate]\label{thm:annular}
Let \(n\ge 2\) and \(p_c=2n/(n-1)\). For every \(R\ge 2\),
\[
  \big\|e^{-it\sqrt{-\Delta}}\phi_R f\big\|_{L^{p_c}(\mathbb R^n\times[1,2])}
  \lesssim_{n,\phi} [\log(2+R)]^{2/p_c}
  \|f\|_{L^{p_c}(\mathbb R^n)}.
\]
\end{theorem}

The plus sign is again covered by complex conjugation, which replaces \(\phi(\eta)\) by \(\overline{\phi(-\eta)}\), a cutoff of the same class. Some logarithmic loss is necessary: for suitable \(\phi\), the focusing example of~\cite[Proposition~8.1]{openai}, recalled at the end of Section~\ref{sec:annular-proof}, shows that the constant in Theorem~\ref{thm:annular} grows at least like \((\log R)^{1/p_c}\).

OpenAI's proof~\cite[Sections~2--7 and Appendix~A]{openai} hinges on Theorem~\ref{thm:annular}. The argument first takes place on a fixed spatial torus, where one frequency annulus contains finitely many modes. It then returns to \(\mathbb R^n\), interpolates with the \(L^2\) and fixed-time bounds, and sums over annuli to prove Theorem~\ref{thm:local-smoothing}. Section~\ref{sec:preliminaries} develops the tools: a conserved current built from the wave, the quaternionic Hessian with its determinant and cofactor, and a Dirichlet problem. Section~\ref{sec:proof} combines them to prove both theorems.

The proof of the annular estimate has three parts, which we now describe informally.

\emph{A conserved positive tensor.} An imaginary displacement of a wave in a light-ray direction puts a Gaussian weight on the contributing frequencies, and so acts as a filter that selects a frequency cap. The displacements are chosen through a symmetric first-order form of the wave equation and depend on auxiliary variables, in the end two real ones for each of the \(n+1\) spacetime coordinates. The outcome is a positive semidefinite matrix field on the enlarged space, the \emph{current}, whose columns are divergence free in the physical and auxiliary variables together. Its \(tt\) entry, the diagonal entry in the time coordinate, measures the square of the filtered wave. The square is taken after the filtered waves are added, so the current still detects their interference. In this sense the energy that we build from the current below is \emph{coherent}.

\emph{A determinant equation on fibers.} Let \(\tau\in\mathbb R^2\) be the pair of auxiliary variables attached to \(t\), and call \(z=(t,\tau)\) the \emph{base point}. The \emph{fiber} over \(z\) consists of the remaining variables: the spatial point \(x\) and the \(2n\) real auxiliary variables attached to the spatial coordinates. On a domain in each fiber we solve a Monge--Ampère equation whose right-hand side is a constant multiple of the \emph{density} \(\mu\), the \(tt\) entry of the current, with zero data on the boundary of the domain, the \emph{fiber wall}. The Hessian in this equation is the quaternionic Hessian of Section~\ref{sec:hessian}, whose determinant has degree \(n\) on a fiber. The unknown is a correction \(u\le0\) to a fixed background potential whose Hessian has rank one, and we call \(-u\) the \emph{depth}. The background supplies one of the \(n\) factors of the determinant. The other \(n-1\) factors turn a local lower bound \(\mu\) for the density into a depth of order \(\mu^{1/(n-1)}\), up to a factor that depends on the size of the comparison ball. Multiplying the depth by the density gives \(\mu^{n/(n-1)}\), and the \emph{energy} of a fiber is the integral of this product. Since the density measures a squared wave, this is how \(p_c=2n/(n-1)\) enters.

\emph{An energy inequality and a scalar comparison.} Because the current is conserved in all variables, we can compute the Laplacian \(\Delta_z\) of the energy in the three base variables. A correction to the energy cancels the one term without a sign, the determinant of the full Hessian. What remains is a term with a favorable sign and a flux through the fiber wall, where the Gaussian is small. The corrected energy \(F\) therefore satisfies \(\Delta_zF\le VF\) with \(V=C|\tau|^{-2}e^{-cR|\tau|}\). On circles \(|\tau|\asymp1/R\) the energy is comparable, up to a normalizing factor, to the \(p_c\)-th power of the \(L^{p_c}\) norm of the wave at the given time, and a comparison argument for the scalar inequality bounds the energy at times \(1\le t\le2\) by the energy at time \(0\). The comparison passes from radius \(|\tau|\asymp1/R\) out to radius one and back. Each of the two passages costs one logarithm, because the radial harmonic functions in two dimensions are affine in the logarithm of the radius. These two logarithms appear in the estimate for the \(p_c\)-th power of the norm.

In more familiar terms, at auxiliary radius \(s=|\tau|\) the filter resolves frequency directions at angular scale \((Rs)^{-1/2}\). This is a fixed scale when \(s\asymp1/R\) and the wave packet scale \(R^{-1/2}\) when \(s=1\). Section~\ref{sec:cylinder} compares the passage between these radii with an induction on scales.

One possible reason for the effectiveness of OpenAI's method, which we refer to as the \emph{coherent energy} method, is that it may analyze the interference between waves more precisely than an induction on scales. The current is quadratic in the sum of the filtered waves, so the cross terms between different waves are part of the conserved quantity; see the remark on interference in Section~\ref{sec:current}. The energy balance of Lemma~\ref{lem:energy-balance} is an identity, and the only term discarded from it afterwards has a sign. An induction on scales, by contrast, passes the interference from one scale to the next through a chain of inequalities, and each of them may lose part of the information about cancellation between waves.

\subsection{A real model}\label{sec:real-model}

The mechanism behind Theorem~\ref{thm:annular} has a real precedent, Serre's compensated integrability~\cite{serre-divergence}. Let \(\Omega\subset\mathbb R^n\), \(n\geq2\), be a bounded convex domain, and let \(A=(A_{ab})\) be a smooth field of positive semidefinite symmetric \(n\times n\) matrices on \(\overline\Omega\). Its \emph{column divergence} is the vector with entries
\[
 (\operatorname{Div}A)_b=\sum_{a=1}^n\partial_{x_a} A_{ab},
\]
the divergences of the columns of \(A\). Assume that \(\operatorname{Div}A=0\). Let \(\mathbf n\) be the outward unit normal, which exists at almost every point of \(\partial\Omega\). By the symmetry of \(A\), the entries of the vector \(A\mathbf n\) are the outward fluxes of the columns. We write \(|A\mathbf n|\) for its Euclidean length, \(dS\) for surface measure, \(\mathbb B^n=\{x:|x|<1\}\), and \(|\mathbb B^n|\) for the volume of the ball. Serre proved
\begin{equation}\label{eq:serre}
  \int_\Omega(\det A)^{\frac{1}{n-1}}\,dx
  \le\left(\frac{1}{n|\mathbb B^n|^{1/n}}
       \int_{\partial\Omega}|A\mathbf n|\,dS\right)^{\frac{n}{n-1}}.
\end{equation}
For \(A=I\) this is the sharp isoperimetric inequality for convex domains. The word \emph{compensated} refers to the gain supplied by the differential constraint: positivity alone bounds \((\det A)^{1/n}\) by \(\operatorname{tr}(A)/n\), whereas column conservation yields control of the higher power in~\eqref{eq:serre} by the boundary flux.

The proof transports a density to a ball. Set \(g=(\det A)^{1/(n-1)}\). If \(\int_\Omega g\,dx=0\) the estimate is immediate. Otherwise choose \(\rho>0\) with \(|\mathbb B^n|\rho^n=\int_\Omega g\,dx\), so that the measure \(g\,dx\) on \(\Omega\) has the same mass as Lebesgue measure on the ball \(\rho\mathbb B^n\). Among the maps \(T\) that push the first measure forward to the second, the \emph{Brenier map} minimizes \(\int_\Omega |x-T(x)|^2g(x)\,dx\). Since the source has a density and both measures have bounded support, this optimal map exists, is unique almost everywhere with respect to \(g\,dx\), and equals \(\nabla\psi\) for a convex function \(\psi\)~\cite{brenier}.

We give the calculation for smooth \(\psi\) and positive \(g\), and do not address the passage to degenerate densities or boundary traces. The change of variables formula for a transport to density one gives \(g=\det DT=\det D^2\psi\). This is the \emph{real Monge--Ampère equation}, which prescribes the determinant of the Hessian of a convex function. The target ball bounds the size of the transport, \(|\nabla\psi|\leq\rho\). The arithmetic--geometric mean inequality for the eigenvalues of \(A^{1/2}D^2\psi A^{1/2}\) gives
\[
  g=(g\det A)^{1/n}=\big(\det(AD^2\psi)\big)^{1/n}
   \le\frac1n\operatorname{tr}(AD^2\psi).
\]
Since the columns of \(A\) are divergence free, \(\operatorname{tr}(AD^2\psi)=\operatorname{div}(A\nabla\psi)\). Integrating over \(\Omega\) gives
\[
  n\int_\Omega g\,dx
  \le\int_{\partial\Omega}\nabla\psi\cdot A\mathbf n\,dS
  \le\rho\int_{\partial\Omega}|A\mathbf n|\,dS ,
\]
which is~\eqref{eq:serre} because \(|\mathbb B^n|\rho^n=\int_\Omega g\,dx\).

Section~\ref{sec:proof} uses the same three ingredients: a positive tensor with divergence-free columns, a Monge--Ampère potential whose determinant is read off from the tensor, and an integration by parts whose interior term has a sign and whose boundary term is a flux. The tensor is the current of Section~\ref{sec:current}, and the density \(\mu\) is its \(tt\) entry, after a spatial averaging, a cutoff and a small positive regularization. The potential solves a quaternionic Monge--Ampère equation on each fiber (Section~\ref{sec:fiber-potentials}); no transport map is involved. The integration by parts is the energy balance of Lemma~\ref{lem:energy-balance}, in which a Schur complement isolates the interior term with a sign and the boundary term is the flux through the fiber wall. The exponent also arises in the same way: if \(A\) is comparable to \(\mu I\), the left side of~\eqref{eq:serre} is comparable to \(\int_\Omega\mu^{n/(n-1)}\,dx=\int_\Omega\mu^{p_c/2}\,dx\).

The quaternionic argument differs from the real one in two main ways. First, its conservation law holds in all variables, while the potential lives on fibers. Integration by parts on a fiber therefore leaves the Laplacian of the energy in the three base variables, and the estimate becomes a differential inequality, which Section~\ref{sec:cylinder} integrates by a scalar comparison. Second, the pairing of the current with the full Hessian contains the full determinant, which has no sign. The correction to the energy in Section~\ref{sec:energy} removes it.

\subsection{A reading guide}\label{sec:reading-guide}

The core of the argument consists of Lemma~\ref{lem:energy-balance}, the energy balance, and Proposition~\ref{prop:fiber-estimates}, the fiber estimates. A remark in Section~\ref{sec:energy} carries out the energy balance in a model with two real variables. On a first reading one can take the following for granted.
\begin{itemize}
\item Proposition~\ref{prop:wave-current} and the formula~\eqref{eq:time-component} for the \(tt\) entry of the current. The Gram row calculation in Section~\ref{sec:conservation} and the proof of the proposition can be skipped. Section~\ref{sec:euclidean} explains where the Euclidean structure enters.
\item The quaternionic linear algebra: Lemmas~\ref{lem:hessian}--\ref{lem:cofactor-divergence} and the boundary identities~\eqref{eq:green} and~\eqref{eq:boundary-compression}.
\item The Dirichlet problem. It enters only through the properties listed after Theorem~\ref{thm:dirichlet}, chiefly the comparison principle of Lemma~\ref{lem:comparison}. The proof of the theorem is in Appendix~\ref{app:dirichlet}.
\end{itemize}
Lemma~\ref{lem:cylinder} is a comparison argument for a scalar elliptic inequality, and Section~\ref{sec:cylinder} relates it to induction on scales. The passage from the torus to \(\mathbb R^n\) in Section~\ref{sec:annular-proof} and the interpolation in Section~\ref{sec:interpolation} are standard. The focusing example at the end of Section~\ref{sec:annular-proof} tests the logarithms in Theorem~\ref{thm:annular}.

\subsection{AI disclosure}
ChatGPT and Claude were used in preparing this note, at the author's direction. They expanded the explanations of material that harmonic analysts may find unfamiliar, related parts of the argument to induction on scales and wave packets, edited the exposition, and added figures.

\section{Preliminaries}\label{sec:preliminaries}

Write \(d=n-1\). Table~\ref{tab:notation} lists the main notation. Throughout, the letters \(i,j,k\) are reserved for the quaternion units, coordinates are indexed by \(a,b\), and the real slots of a coordinate, introduced in Section~\ref{sec:current}, are indexed by \(\alpha,\beta\).

\begin{table}[!htbp]
\centering
\caption{Main notation.}
\label{tab:notation}
\small
\renewcommand{\arraystretch}{1.15}
\begin{tabular}{@{}p{.31\linewidth}p{.65\linewidth}@{}}
\hline
Symbol & Meaning \\
\hline
\(R,\phi_R,p_c,\sigma_p\) & Annular frequency, annular multiplier, critical exponent and Sobolev order in Theorem~\ref{thm:local-smoothing}. \\
\(\bar t,\bar x\) & Real null coordinates \((x_n\pm t)/\sqrt2\). \\
\(\kappa,v\) & Cone coordinates of a frequency \(\xi\): \(\kappa=(\xi_n+|\xi|)/\sqrt2\) and the slope \(v=\xi'/\kappa\). \\
\(w_\ell,\widetilde f,\mathsf q_I\) & The two states (\(\ell=1,2\)), translated system solution and quaternionic Gram row. \\
\(\vartheta,\mathcal J_\Theta\) & Unit phase \(r_0/|r_0|\) and absolute real Jacobian of \(\Theta\). \\
\(\tau,P,Z\) & Auxiliary variables attached to \(t\), \(x'\) and \(x_n\). \\
\(Y,Z_-,y\) & Auxiliary null coordinates \((Z\pm\tau)/\sqrt2\), and \(y=|Y|\). \\
\(W,\Phi\) & Angular coordinate \(P/Y\) and passive coordinate. \\
\(A_0,A_w,A\) & Wave current, its smoothed and cut-off version, and the regularized current. \\
\(G,H,C\) & Encoded gradient, quaternionic Hessian and Moore cofactor. \\
\(g\) & Prescribed determinant in the Dirichlet problem. \\
\(z,s\) & Base point \((t,\tau)\) and auxiliary radius \(|\tau|\). \\
\(\gamma\) & Aperture of the frequency cap \(|v|\le\gamma\). \\
\(K,a,Q,\chi_0\) & Large constant, fiber scale \(a=K/R\), normalized coordinate \(\Phi/a\) and its cutoff. \\
\(r,\Omega_z\) & Defining function and the fiber domain \(\Omega_z=\{r<0\}\) over \(z\). \\
\(\mathcal K_\delta\) & Periodized averaging kernel at spatial scale \(\delta\). \\
\(\lambda,\mu,B_0,u\) & The constant \(R^{30n}\), density \(A_{tt}\), rank-one background potential \(\lambda|\Phi|^2\) and Dirichlet correction. \\
\(h,\Sigma\) & \(tt\) entry of the full Hessian and positive semidefinite remainder in the Schur decomposition. \\
\(F_0,E,F,\overline F\) & Raw energy, its correction, the energy \(F=F_0+E\) and its circular mean. \\
\(M,\Pi\) & Averaging constant and averaged profile in the proof of Proposition~\ref{prop:fiber-estimates}. \\
\(b\) & Inner radius of the cylinder in Lemma~\ref{lem:cylinder}; \(b=2a\) in the application. \\
\(c_0,C_0,c_{\mathrm{wall}}\) & Annular support constants and Gaussian decay rate at the wall. \\
\hline
\end{tabular}
\end{table}

\subsection{A positive current that keeps the interference between waves}\label{sec:current}

In this note a \emph{current} is a matrix field that represents a conservation law, column by column: in real coordinates \(X\) the equation is \(\sum_I\partial_{X_I} A^{\mathbb R}_{IJ}=0\) for each \(J\), where the derivatives run over the auxiliary variables as well as physical spacetime. (Currents in the sense of de~Rham, that is, distributions acting on differential forms, play no role here.)

The construction starts from the quadratic flux of a symmetric hyperbolic
system, a special case of the bilinear flux of Friedrichs~\cite{friedrichs}.
Kaiser and Streater~\cite{kaiser-streater} used imaginary spacetime
translations to analyze waves, and Kaiser~\cite[Section~4.5, (36)--(40)]{kaiser}
constructed a conserved current in physical spacetime for positive-energy
Klein--Gordon waves. OpenAI~\cite[Section~2]{openai} uses the quadratic flux
to select a Gaussian filter and then obtains a positive tensor with
divergence-free columns in the combined physical and auxiliary variables.
The auxiliary variables used in the filter will also define the fibers.

The construction has three steps. First, the matrix symbol of the wave
equation provides a Gaussian cutoff in slope. Second, we build a positive real
matrix whose columns have zero divergence in the physical and auxiliary
variables together. Third, we change auxiliary coordinates so that each
physical coordinate, together with two real auxiliary coordinates, forms a
single quaternionic coordinate. We call these three real coordinates the
\emph{slots} of the quaternionic coordinate. In the new coordinates the
\(tt\) entry of the matrix has a direct description as a Fourier multiplier,
and the other entries will allow us to integrate that entry against a
potential on a fiber.

\subsubsection{The symbol and its Gaussian.}\label{sec:symbol}
Use an orthonormal spatial frame \(x=(x',x_n)\) with \(x'\in\mathbb R^d\).
In the null coordinates \(\bar{t}=(x_n+t)/\sqrt2\) and
\(\bar{x}=(x_n-t)/\sqrt2\), the wave operator \(\Delta_x-\partial_t^2\) becomes
\(2\partial_{\bar{t}}\partial_{\bar{x}}+\Delta_{x'}\). Indeed,
\(\partial_{x_n}=(\partial_{\bar t}+\partial_{\bar x})/\sqrt2\)
and \(\partial_t=(\partial_{\bar t}-\partial_{\bar x})/\sqrt2\), so
the diagonal second derivatives in these two squares cancel.

The cone has a useful parametrization on the patch \(\xi_n+|\xi|>0\).
The coefficients of \(\bar x\) and \(\bar t\) in \(x\cdot\xi-t|\xi|\)
are \((\xi_n+|\xi|)/\sqrt2\) and \((\xi_n-|\xi|)/\sqrt2\).
Multiplying them gives \(-|\xi'|^2/2\) by the cone equation
\(|\xi|^2=\xi_n^2+|\xi'|^2\). We can therefore write a mode as
\[
  u_{\kappa,v}
    =c_{\kappa,v}\exp\!\left(i\kappa
       (\bar{x}+v\cdot x'-\tfrac12|v|^2\bar{t})\right),
  \qquad
  \kappa=\frac{\xi_n+|\xi|}{\sqrt2},
  \quad
  v=\frac{\xi'}{\kappa}.
\]
Near the axis \(\xi'=0\), \(v\) is a slope coordinate on the cone, and
differences of slopes are comparable to angular distance.

We now put the scalar wave into a first-order system. Its state has \(n\)
components, indexed by \(\rho,\sigma\in\{0,\dots,d\}\), and \(E_{\rho\sigma}\)
denotes the elementary matrices. Let \(\partial_a\) denote the derivative in
the \(a\)-th coordinate of the null frame \((\bar{t},x',\bar{x})\). The
symmetric system \(\mathcal P=\sum_a\mathsf C^a\partial_a\) has coefficients
\[
  \mathsf C^{\bar{t}}=2E_{00},\qquad
  \mathsf C^{x_\rho}=E_{0\rho}+E_{\rho0}\quad(1\le\rho\le d),\qquad
  \mathsf C^{\bar{x}}=-\sum_{\rho=1}^d E_{\rho\rho}.
\]
Applied to \((\partial_{\bar{x}}\psi,\nabla_{x'}\psi)^{\mathsf t}\),
the first row is \(2\partial_{\bar t}\partial_{\bar x}\psi+\Delta_{x'}\psi\).
Row \(\rho\) for \(1\leq \rho\leq d\) is
\(\partial_{x_\rho}\partial_{\bar x}\psi-\partial_{\bar x}\partial_{x_\rho}\psi\),
which is zero. Thus a scalar wave gives a solution of this system.

For a mode of covector \(p=\kappa(-|v|^2/2,v,1)\), the matrix symbol with
the common factor \(i\) removed is
\[
  \sum_a p_a\mathsf C^a
     =\kappa\begin{pmatrix}-|v|^2&v^{\mathsf t}\\v&-I\end{pmatrix}.
\]
Its negative quadratic form is a square:
\[
  -w^{\mathsf t}\!\left(\sum_a p_a\mathsf C^a\right)w
     =\kappa|w'-vw_0|^2
  \qquad (w=(w_0,w')\in\mathbb R^n).
\]
(Expanded, the left side is
\(\kappa(|v|^2w_0^2-2w_0\,v\cdot w'+|w'|^2)\).)
In particular \((1,v)^{\mathsf t}\) spans the symbol kernel, so
\(
  f_{\mathrm{sys}}
     =\sum_{\kappa,v}\kappa^{d/2}(1,v)^{\mathsf t}u_{\kappa,v}
\)
is an entire solution of \(\mathcal P f_{\mathrm{sys}}=0\) for any finite
sum of modes. The vector \((1,v)^{\mathsf t}\), called the polarization,
records the relative amplitudes of the entries of the first-order system.
Every entry carries the same scalar oscillation.

Take two real states \(w_1,w_2\) and translate the physical coordinates
by \(-i\mathbf q\), where
\[
  \mathbf q^a
    =\frac12\sum_{\ell=1}^2 w_\ell^{\mathsf t}\mathsf C^a w_\ell .
\]
Substituting the complex coordinates into a mode adds \(p\cdot\mathbf q\)
to the exponent, because \(i(-i)=1\). The square identity above gives
\[
  p\cdot\mathbf q=-\frac{\kappa}{2}
       \sum_{\ell=1}^2 |w_\ell'-v w_{\ell0}|^2 .
\]

This displacement has a simple geometric meaning. The numbers
\(w^{\mathsf t}\mathsf C^aw\) are the values on a constant state \(w\) of
the quadratic flux \(f\mapsto(f^{\mathsf t}\mathsf C^af)_a\) of the system,
so \(\mathbf q\) is half the sum of the fluxes of the two states. For a
single state \(w=(w_0,w')\), half the flux is the vector with null-frame
components
\[
  \tfrac12\bigl(w^{\mathsf t}\mathsf C^{\bar t}w,\,
     (w^{\mathsf t}\mathsf C^{x_\rho}w)_\rho,\,
     w^{\mathsf t}\mathsf C^{\bar x}w\bigr)
  =\bigl(w_0^2,\,w_0w',\,-\tfrac12|w'|^2\bigr).
\]
Since \(x_n^2-t^2=2\bar t\bar x\), the Minkowski form \(|x|^2-t^2\) equals
\(|x'|^2+2\bar t\bar x\) in the null frame, and it vanishes on this vector.
Each state therefore contributes an imaginary displacement along a light
ray. For \(w_0\ne0\) it is the ray of the mode with slope \(v=w'/w_0\),
because the ray \((t,x)=(1,\xi/|\xi|)\) of a mode has null-frame components
proportional to \((1,v,-\tfrac12|v|^2)\). By the square identity this mode
is not damped at all, while a mode of slope \(v\) is multiplied by
\(\exp(-\tfrac12\kappa w_0^2|v-w'/w_0|^2)\), a Gaussian in the difference
of the two slopes. With two states the two weights multiply.

To write the sum of the two squares as a single complex norm we pair the
two states, and to keep this pairing apart from the complex structure of
the Fourier phase we do so inside the quaternions. Recall that a
quaternion is a real linear combination of \(1,i,j,k\), where
\(i^2=j^2=k^2=-1\), \(ij=k\), \(jk=i\), \(ki=j\), and reversing any of
these last three products changes its sign. The two planes
\(\mathbb C_i=\mathbb R+i\mathbb R\) and \(\mathbb C_k=\mathbb R+k\mathbb R\)
are each copies of the complex numbers, but \(ik=-ki\). We identify \(i\)
with the imaginary unit of the Fourier phase and pair the two states as
\(w_1+kw_2=(r_0,\zeta)\in\mathbb C_k\times\mathbb C_k^d\).
The multiplier on each mode is then
\[
  e^{p\cdot\mathbf q}
       =\exp\!\left(-\tfrac12\kappa |\zeta-vr_0|^2\right).
\]
For fixed \(r_0\ne0\) this filter is centered at the slope \(\zeta/r_0\) and
selects the modes whose slope lies within about \((\kappa|r_0|^2)^{-1/2}\)
of it. It is a cap cutoff whose width varies with the auxiliary coordinate.
For a fixed mode and fixed \(r_0\), the squared weight in \(\zeta\) has integral
\[
  \int_{\mathbb C_k^d}\kappa^d e^{-\kappa|\zeta-vr_0|^2}\,d\zeta
    =\int_{\mathbb C_k^d}e^{-|\zeta|^2}\,d\zeta ,
\]
by the substitution \(\zeta\mapsto\sqrt\kappa(\zeta-vr_0)\). This explains the
normalizing power \(\kappa^{d/2}\) in the amplitude.

\subsubsection{Conservation before and after the change of variables.}\label{sec:conservation}
The Gaussian localization alone gives a family of filtered waves.
To use them in the proof, we need a differential identity that involves
both the Gaussian parameters and the physical variables.
Set \(\widetilde f=f_{\mathrm{sys}}((\bar{t},x',\bar{x})-i\mathbf q)\) and define the vector fields
\[
  D_{\ell\rho}
    =\partial_{w_{\ell\rho}}
      +i\sum_a(\mathsf C^a w_\ell)_\rho\,\partial_a,
  \qquad
  L_\ell=\sum_\rho\overline{\widetilde f_\rho}\,D_{\ell\rho}.
\]
The \(D_{\ell\rho}\) have the form of the complex vector fields associated with
tube coordinates in Baouendi and Trèves~\cite[Section~1]{baouendi-treves};
everything we need about them is computed below.
By the symmetry of \(\mathsf C^a\),
\(\partial_{w_{\ell\rho}}\mathbf q^a=(\mathsf C^a w_\ell)_\rho\).
The two terms in \(D_{\ell\rho}\) therefore cancel when applied to the
\(a\)-th component of \((\bar t,x',\bar x)-i\mathbf q\).
Since \(f_{\mathrm{sys}}\) is holomorphic in these translated variables,
\(D_{\ell\rho}\widetilde f=0\). Conjugating gives
\(\overline{D_{\ell\rho}}\,\overline{\widetilde f}=0\), so
\[
  D_{\ell\rho}\overline{\widetilde f}
     =(D_{\ell\rho}-\overline{D_{\ell\rho}})\overline{\widetilde f}
     = 2 i\sum_a (\mathsf C^a w_\ell)_\rho\,
             \partial_a\overline{\widetilde f}.
\]
The state coefficients of \(D_{\ell\rho}\) are constant and its
physical coefficients are independent of the physical variables,
so its divergence is zero. By the product rule,
\[
  \operatorname{div}L_\ell
   =\sum_\rho D_{\ell\rho}\overline{\widetilde f_\rho}
   =2i w_\ell^{\mathsf t}\mathcal P\overline{\widetilde f}=0 .
\]
For the last equality, note that for fixed states the translation is
constant in the physical variables and so commutes with \(\mathcal P\),
and that the coefficients of \(\mathcal P\) are real, so the equation
\(\mathcal P\widetilde f=0\) can be conjugated.

We will also need the derivative of \(\overline{L_{\ell'}}\) along \(L_\ell\).
Here \(\nabla_{L_\ell}\overline{L_{\ell'}}\) means that \(L_\ell\)
differentiates the coefficients of the vector field
\(\overline{L_{\ell'}}=\sum_\sigma\widetilde f_\sigma\overline{D_{\ell'\sigma}}\).
The derivative of the amplitude \(\widetilde f_\sigma\) vanishes. In the
remaining term, \(D_{\ell\rho}\) acts on the coefficient
\(-i(\mathsf C^a w_{\ell'})_\sigma\) of \(\overline{D_{\ell'\sigma}}\),
giving \(-i\delta_{\ell\ell'}\mathsf C^a_{\sigma\rho}\). Consequently
\[
  \nabla_{L_\ell}\overline{L_{\ell'}}
     =-i\delta_{\ell\ell'}\sum_a
            (\overline{\widetilde f}^{\mathsf t}\mathsf C^a \widetilde f)\,\partial_a.
\]
Each quadratic form in parentheses is real because \(\mathsf C^a\) is
real symmetric; it is the \(a\)-th component of the quadratic flux of the
translated solution \(\widetilde f\). For a real or complex matrix field \(B\), the column
divergence has \(J\)-th entry \(\sum_I\partial_I B_{IJ}\), the ordinary
divergence of the vector field given by column \(J\).
The column divergence of an outer product is therefore
\[
  \operatorname{div}(L_\ell\otimes\overline{L_\ell})
    =(\operatorname{div}L_\ell)\overline{L_\ell}
        +\nabla_{L_\ell}\overline{L_\ell}.
\]
The first term vanishes and the second is purely imaginary. It follows
that \(\mathsf B=\operatorname{Re}\sum_{\ell=1}^2
L_\ell\otimes\overline{L_\ell}\) has divergence-free real columns.
It is also positive semidefinite: writing a complex field as its real
part plus \(i\) times its imaginary part shows that the real part of its
outer product with its conjugate is the sum of two real outer squares,
so \(\mathsf B\) is a sum of four real outer squares.

At this point we have a positive, conserved tensor in the state variables.
The rest of this subsection checks two things: that in the new coordinates
the four squares form a single quaternionic outer product, and that the
change of variables introduces no divergence error. The outcome is
Proposition~\ref{prop:wave-current}, and the check can be skipped on a
first reading.

To see how the four real squares fit together, write
\(L_\ell=\sum_I L_\ell^I\partial_I\), where \(I\) runs over all the real
coordinates, and set \(\mathsf q_I=L_1^I+kL_2^I\).
The coefficients \(L_\ell^I\) lie in \(\mathbb C_i\).
The planes \(\mathbb C_i\) and \(k\mathbb C_i\) are orthogonal for
the real inner product \(\operatorname{Re}(q\overline{q'})\) on quaternions:
the first plane is spanned by \(1,i\), and the
second by \(k,j\), since \(ki=j\).
Left multiplication by \(k\) is an isometry. In the inner product
of two of the quaternions \(\mathsf q_I\), the mixed terms between the
two planes therefore vanish, so
\[
  \operatorname{Re}(\mathsf q_I\overline{\mathsf q_{I'}})
   =\operatorname{Re}\!\left(
      L_1^I\overline{L_1^{I'}}+L_2^I\overline{L_2^{I'}}\right)
   =\mathsf B_{II'}.
\]
Thus \(\mathsf B\) is the Gram matrix of the quaternions \(\mathsf q_I\),
and we call \(\mathsf q_I\) the \emph{Gram row} of the coordinate \(I\).
The Gram rows of the state coordinates are
\[
  (\mathsf q_{w_{1\rho}},\mathsf q_{w_{2\rho}})
   =(\overline{\widetilde f_\rho},k\overline{\widetilde f_\rho}).
\]

Adjoin a passive coordinate \(\eta\in\mathbb C_k\) and set
\[
  \omega^a=\tfrac12 w^{\mathsf t}\mathsf C^a w
                +\delta_{a,\bar{x}}\eta,
  \qquad w=w_1+kw_2.
\]
Write \(\Theta\) for this map together with the identity on the real
spacetime variables. The Gram rows in the \(\eta\) directions are zero. The
state derivatives satisfy \(\partial_{w_\rho}\omega^a=(\mathsf C^a w)_\rho\),
so the Gram row of each real spacetime coordinate is
\[
  \mathsf q_a=\sum_\rho (\partial_{w_\rho}\omega^a)\,i\,\overline{\widetilde f_\rho}.
\]
To compute the Gram rows of the two new auxiliary coordinates over a
spacetime coordinate, recall that under a change of variables the
components of a vector field transform by the differential. The same is
true of the Gram rows: the row of a new real coordinate \(\Theta_J\) is
\(\sum_I(\partial_I\Theta_J)\,\mathsf q_I\). Write
\(\omega^a=\omega_1^a+k\omega_2^a\), so that
\(\omega^a i=\omega_1^a i+\omega_2^a j\).
Since \(\omega^a\) is a holomorphic polynomial in the variables
\(w_\rho=w_{1\rho}+kw_{2\rho}\), its real differential in the
\((w_{1\rho},w_{2\rho})\) variables is multiplication by
\(c=\partial_{w_\rho}\omega^a=c_1+kc_2\), whose real matrix appears on the
left below. With \(q=\overline{\widetilde f_\rho}\), applying this matrix
to the two Gram rows in that state pair gives
\[
  \begin{pmatrix}c_1&-c_2\\c_2&c_1\end{pmatrix}
  \begin{pmatrix}q\\kq\end{pmatrix}
  =\begin{pmatrix}\overline c q\\k\overline c q\end{pmatrix}.
\]
The first of these two entries equals \(-i(ciq)\), since
\(-ici=\overline c\) for \(c\in\mathbb C_k\).
The second entry equals \(-j(ciq)\), since \(ki=j\).
Summing over \(\rho\) thus gives the three rows
\(\mathsf q_a,-i\mathsf q_a,-j\mathsf q_a\) over the real
coordinate \(a\) and its two auxiliary slots.
In other words, with \((e_0,e_1,e_2)=(1,i,j)\), the Gram row for slot
\(\alpha\) of the coordinate \(a\) is \(\overline{e_\alpha}\mathsf q_a\).
In the real Gram product of two such rows, quaternionic conjugation
reverses the order of the factors, so the Gram entry for the slots
\((a,\alpha)\) and \((b,\beta)\) is
\(\operatorname{Re}(\overline{e_\alpha}\,\mathsf q_a\,
\overline{\mathsf q_b}\,e_\beta)\).
In general, the \emph{real representation} in three slots of a
quaternionic matrix \(A=(A_{ab})\) is the real matrix
\begin{equation}\label{eq:real-representation}
  A^{\mathbb R}_{a\alpha,b\beta}
     =\operatorname{Re}(\overline{e_\alpha}A_{ab}e_\beta).
\end{equation}
The Gram entries just computed are the real representation of
the single quaternionic outer product
\(\big(\mathsf q_a\overline{\mathsf q_b}\big)_{a,b}\).

On an inverse chart, the change of variables sends a
vector to its image under \(D\Theta\) and an outer product
of two vectors to the outer product of their images. We also carry
along the source volume element. The transformed matrix field
relative to the target volume is therefore
\[
  \widetilde{\mathsf B}\circ\Theta
    =\mathcal J_\Theta^{-1}D\Theta\,\mathsf B\,D\Theta^{\mathsf t},
  \qquad \mathcal J_\Theta=|\det\nolimits_{\mathbb R}D\Theta|.
\]
This is the pushforward of \(\mathsf B\) as a tensor density. The factor
\(\mathcal J_\Theta^{-1}\) compensates for the change of volume in
integration. The divergence still acquires an additional term,
because a column in the new coordinates is a position-dependent
combination of the old columns.

To see this in the weak formulation, let \(\varphi\) be a test function
supported in the target chart and fix a real target component \(\Theta_J\).
The integral on the left uses target coordinates, and those on the right
use source coordinates. The chain rule and the change of variables give
\[
 \int (\operatorname{div}\widetilde{\mathsf B})_J\,\varphi\,d\mathrm{vol}
 =-\int\sum_{I,I'}\mathsf B_{II'}\,
       \partial_I(\varphi\circ\Theta)\,\partial_{I'}\Theta_J\,d\mathrm{vol}.
\]
Integrating by parts on the right gives the two terms
\[
 \begin{aligned}
 &\int (\varphi\circ\Theta)\sum_{I'}(\operatorname{div}\mathsf B)_{I'}
      \,\partial_{I'}\Theta_J\,d\mathrm{vol}\\
 &\qquad+\int(\varphi\circ\Theta)\sum_{I,I'}\mathsf B_{II'}\,
          \partial_I\partial_{I'}\Theta_J\,d\mathrm{vol}.
 \end{aligned}
\]
The first integral vanishes by conservation in the source variables.
In the second, only derivatives in the state directions contribute,
since \(\Theta\) is the identity in the physical variables and affine
in the passive variable. The two state blocks of \(\mathsf B\) are both
\(\left(\operatorname{Re}(\widetilde f_\rho\overline{\widetilde f_\sigma})\right)_{\rho\sigma}\)
and their mixed block is zero, so the two block contractions
sum to zero. Indeed, for a holomorphic function of \(w=w_1+kw_2\),
a derivative in \(w_{2\rho}\) is \(k\) times the derivative in
\(w_{1\rho}\). Differentiating twice replaces \(k^2\) by \(-1\), so
on every real component of \(\Theta\)
\[
  \big(\partial_{w_{1\rho}}\partial_{w_{1\sigma}}
      +\partial_{w_{2\rho}}\partial_{w_{2\sigma}}\big)\Theta_J=0.
\]
The real tensor therefore still has divergence-free columns after the change
of variables.

\subsubsection{The current in physical coordinates.}
We now make the change of variables explicit. In the null frame the
quadratic expressions for \(\omega^{\bar t}\), \(\omega^{x_\rho}\) and
\(\omega^{\bar x}\) are respectively \(r_0^2\), \(r_0\zeta_\rho\) and
\(\eta-\tfrac12\zeta\cdot\zeta\), so the auxiliary map is
\[
  (r_0,\zeta,\eta)\longmapsto
  (Y,P,Z_-)=(r_0^2,r_0\zeta,\eta-\tfrac12\zeta\cdot\zeta).
\]
Use \((\tau,P,Z)\) for the auxiliary variables in the physical frame
\((t,x',x_n)\). The resulting single-valued coordinates on \(Y\ne0\) are
\begin{equation}\label{eq:wave-cap-coordinates}
  Y=\frac{\tau+Z}{\sqrt2},\qquad W=\frac{P}{Y},\qquad
  \Phi=\frac{Z-\tau}{\sqrt2}+\frac{P\cdot P}{2Y},\qquad y=|Y|.
\end{equation}
Here \(Z_-=(Z-\tau)/\sqrt2\), and the dot product of two
\(\mathbb C_k\)-vectors in a holomorphic formula is bilinear. The two
local inverses have \(r_0=\pm\sqrt Y\), \(\zeta=r_0W\), and \(\eta=\Phi\).
The complex Jacobian matrix in the displayed order is lower triangular
by blocks, with diagonal blocks \(2r_0\), \(r_0 I\) and \(1\), so
its determinant is \(2r_0^{d+1}\). Each complex coordinate has two
real components, so the absolute real determinant is the modulus squared
and \(\mathcal J_\Theta=4|r_0|^{2d+2}\).
Passing from one inverse to the other replaces \(w\) by \(-w\) and fixes
\(\eta\). The translated wave \(\widetilde f\) stays the same because
\(\mathbf q\) is quadratic in \(w\). The derivatives
\(\partial_{w_\rho}\omega^a\) change sign, so all three transformed rows
over each \(a\) change sign, while their real Gram tensor and the absolute
Jacobian are unchanged. The chartwise densities therefore glue to a single
tensor on \(Y\ne0\).

The state derivative applied to a polarization is
\begin{equation}\label{eq:wave-polarization}
  \begin{gathered}
    (\partial_w\omega)(1,v)^{\mathsf t} = 2r_0 V_v,\qquad
    V_v=\left(1,\frac{v+W}{2},-\frac{W\cdot v}{2}\right)^{\mathsf t},\\
    (\partial\Phi)V_v=0,\qquad (\partial W)V_v=\frac{v-W}{2Y}.
  \end{gathered}
\end{equation}
To check the first formula, vary \((r_0,\zeta)\) in the direction
\((1,v)\) with \(\eta\) fixed. The derivatives of the three components
\((Y,P,Z_-)\) are \(2r_0\), \(\zeta+r_0v\) and \(-\zeta\cdot v\).
Dividing by \(2r_0\) and using \(\zeta=r_0W\) gives \(V_v\).
Since \(\Phi=\eta\) stays fixed, the second formula follows; it expresses tangency to the level sets of \(\Phi\).
The quotient rule gives
\[
  (\partial W)V_v
    =\frac{(v+W)/2-W}{Y}
    =\frac{v-W}{2Y}.
\]
The change in slope along a mode is thus proportional to the mismatch
between the slope of the mode and the center of the filter.

For \(e=Y/y\) define the isometry
\(\iota_e(c_1+ic_2)=c_1+eic_2\) on \(\mathbb C_i\), where \(c_1,c_2\) are real.
For a quaternionic matrix, \({}^*\) denotes the conjugate transpose, and a
matrix \(A\) is Hermitian if \(A=A^*\). The outer product
\(\chi\chi^*\) of a quaternionic column with itself is Hermitian and
positive semidefinite, since its quadratic
form on a column \(v\) is \(|\chi^*v|^2\). Recall
from~\eqref{eq:real-representation} the real representation
\(A^{\mathbb R}\) of \(A\) in three slots.

\begin{proposition}[Wave current]\label{prop:wave-current}
There is a positive semidefinite current \(A_0=\chi\chi^*\), with
\[
  \chi=\sum_{\kappa,v}V_v\,\iota_e(\overline{u_{\kappa,v}})
         \left(\frac{\kappa}{y}\right)^{d/2}
         e^{-\kappa y|W-v|^2/2}.
\]
It is smooth on \(Y\ne0\), \((\partial\Phi)A_0=0\), and the columns of
\(A_0^{\mathbb R}\) have zero real divergence. These properties hold
in both the null and physical frames.
\end{proposition}

In real terms the proposition says the following. The enlarged space has
\(3(n+1)\) real coordinates: the spacetime coordinates \((t,x)\) and the two
real components of each of the \(n+1\) auxiliary variables \((\tau,P,Z)\).
On the open set \(Y\ne0\), \(A_0^{\mathbb R}\) is a smooth field of positive
semidefinite real symmetric matrices of size \(3(n+1)\). Each of its columns
is a divergence-free vector field in the \(3(n+1)\) variables, and each
column is tangent to the level sets of the two real components of \(\Phi\).
The second property will allow us to multiply the current by a cutoff in
\(\Phi\) without losing the first. In the column of the physical time
variable \(t\), the conservation law expresses the \(t\)-derivative of the
\(tt\) entry as minus the divergence, in the other \(3n+2\) variables, of
the remaining entries of that column.

\begin{proof}
By the Gram row calculation, in an inverse chart the transformed real tensor
is the real representation of the quaternionic outer square of
\(\mathcal J_\Theta^{-1/2}(\mathsf q_a)_a\).
For each mode its column is \(V_v \vartheta i\overline{u_{\kappa,v}}\) times
the positive factor in the proposition, where \(\vartheta=r_0/|r_0|\).
This uses \(|\zeta-vr_0|^2=y|W-v|^2\) and
\(2|r_0|\mathcal J_\Theta^{-1/2}=y^{-d/2}\).
Right multiplication of every entry of the column by the same unit
quaternion \((\vartheta i)^{-1}\) preserves its outer square
and changes the scalar factor to
\[
  \vartheta i\overline{u_{\kappa,v}}(\vartheta i)^{-1}
    =\vartheta\overline{u_{\kappa,v}}\,\vartheta^{-1}
    =\iota_e(\overline{u_{\kappa,v}}).
\]
For the first equality, \(i\) commutes with the Fourier amplitude.
For the second, \(i\vartheta^{-1}=\vartheta i\) implies
\(\vartheta i\vartheta^{-1}=\vartheta^2i=ei\).
The resulting formula depends only on \(Y\) and the other target
coordinates, so it is independent of the choice of local square root.

The tangency property \((\partial\Phi)A_0=0\) follows mode by mode
from \eqref{eq:wave-polarization}.
In real variables it says that \(A_0^{\mathbb R}\) kills the
gradients of both real components of \(\Phi\). To see this directly, note
that in the source variables the real and imaginary parts of \(L_\ell\) have zero
\(\eta\) components. Their images under \(D\Theta\) are tangent to
the level sets of \(\Phi\) because \(\Phi\circ\Theta=\eta\). The transformed
real matrix is a sum of the outer squares of these images with a
positive Jacobian factor, so it kills both gradients.

Finally, the change from null to physical
coordinates is a real orthogonal transformation acting on the coordinate
index in all three slots. It preserves both positivity and column
divergence.
\end{proof}

In the estimates we use the fixed spatial torus
\(\mathbb T^n=\mathbb R^n/(100\mathbb Z)^n\). Its Fourier frequencies
are \(\xi\in(2\pi/100)\mathbb Z^n\), so a bounded annulus contains
only finitely many modes. To use the cone parametrization on a cap,
choose an orthonormal frame whose last vector is the axis of the cap.
Then \(\xi_n\) is the component of \(\xi\) along the axis, and \(\xi'\)
consists of its components along the other frame vectors. Expressing
\(x\) in the same frame preserves \(x\cdot\xi\), so every mode remains
periodic with respect to the original lattice \((100\mathbb Z)^n\).
Thus the cap determines the coordinates used to compute \(\kappa\) and
\(v\), while the same lattice indexes the Fourier sum for every cap.

Take the modes of \(e^{-it\sqrt{-\Delta}}f\) in this cap, with
\(c_0R\le |\xi|\le C_0R\), \(|v|\le\gamma\) and \(R\ge2\).
Then \(|\xi|=\kappa(1+|v|^2/2)/\sqrt2\) gives
\(\kappa\asymp R\). There are \(O(R^n)\) frequencies
in the annulus, so Parseval and Cauchy--Schwarz give
\(
  \sum_{\kappa,v}|c_{\kappa,v}|\lesssim
  R^{n/2}\|f\|_{L^2(\mathbb T^n)}.
\)
The current inherits the spatial periods of these modes.

At a fixed auxiliary point, in the coordinates \eqref{eq:wave-cap-coordinates}, take
the physical time component of \(V_v\). Since
\(t=(\bar t-\bar x)/\sqrt2\), this is
\((1+W\cdot v/2)/\sqrt2\). It belongs to \(\mathbb C_k\).
Let \(S_0,S_1\) be the spatial Fourier multipliers obtained by taking
the two real coefficients of this component in the basis \((1,k)\)
and multiplying each by
\((\kappa/y)^{d/2}e^{-\kappa y|W-v|^2/2}\).
By Proposition~\ref{prop:wave-current}, the time component of \(\chi\)
then decomposes as
\[
  \chi_t
    = \iota_e\big(\overline{S_0 e^{-it\sqrt{-\Delta}}f}\big)
       +k\iota_e\big(\overline{S_1 e^{-it\sqrt{-\Delta}}f}\big).
\]
The conjugations here act on the outputs of the operators.
They can be placed outside the sums because the symbols are real.
Conjugation by \(\vartheta\in\mathbb C_k\) carries the two
orthogonal planes \(\mathbb C_i\) and \(k\mathbb C_i\) to
\(\iota_e(\mathbb C_i)\) and \(k\iota_e(\mathbb C_i)\).
It is an isometry, so \(|\chi_t|^2\) is the sum of the squares of the
two scalar outputs. Thus
\begin{equation}\label{eq:time-component}
  A_{0,tt}
    =\big|S_0 e^{-it\sqrt{-\Delta}}f\big|^2
       +\big|S_1 e^{-it\sqrt{-\Delta}}f\big|^2.
\end{equation}
\begin{remark}[Interference between two modes]
Each square in \eqref{eq:time-component} is taken after adding the
filtered modes. Suppose there are just two modes, and let
\(\chi^{(1)},\chi^{(2)}\) be their contributions to the column in
Proposition~\ref{prop:wave-current}.
Since \(A_0=\chi\chi^*\), its \(tt\) entry is
\[
 A_{0,tt}
   = |\chi_t^{(1)}|^2+|\chi_t^{(2)}|^2
     +2\operatorname{Re}\bigl(\chi_t^{(1)}\overline{\chi_t^{(2)}}\bigr).
\]
Changing the sign of the second wave reverses the sign of the last term
while preserving both squares. At a point where the two filtered time
components are equal and opposite, this term cancels both squares exactly.

The interference term itself satisfies a conservation law. Apply
Proposition~\ref{prop:wave-current} to the sum of the two modes and
subtract the identities for the individual modes. The two single-mode currents
drop out, leaving
\[
 \operatorname{div}\!\left[
   \bigl(\chi^{(1)}(\chi^{(2)})^*
        +\chi^{(2)}(\chi^{(1)})^*\bigr)^{\mathbb R}
 \right]=0 .
\]
Here divergence acts on each column in the physical and auxiliary
variables together. In the column corresponding to physical time,
the \(tt\) entry is precisely the signed cross term above. Its time
derivative is balanced by the divergence of the other entries of that
column. Thus the integration by parts used later applies even to the
interference between a single pair of modes.
\end{remark}

The symbol of \(S_0\) contains the factor
\((1+\operatorname{Re}W\cdot v/2)/\sqrt2\), bounded away from zero
for small \(|W|\) and \(|v|\). Where the Gaussian is also bounded
from below, this multiplier can be inverted on the cap to recover the
wave. On parts of the fiber wall the center \(W\) is separated
from the input slopes, and the Gaussian gives exponential decay in
\(Ry\). We will use both features of the same filter.
The other components of \(\chi\) have the same Gaussian with smooth
polarization factors, and the entries of \(A_0\) are their pairwise
products.

\subsubsection{Where the Euclidean structure enters.}\label{sec:euclidean}
As recalled in Section~\ref{sec:earlier}, in odd dimensions \(n\geq3\) the range \(p\geq2+8/(3n-3)\) is sharp for compact Riemannian manifolds, and it lies strictly above \(p_c\). An argument for Theorem~\ref{thm:annular} must therefore use the Euclidean structure somewhere, and readers who know the variable-coefficient counterexamples will want to see where. It enters most visibly through the conservation law of Proposition~\ref{prop:wave-current}. The modes are exact exponentials for the constant-coefficient system \(\mathcal P\). Hence \(f_{\mathrm{sys}}\) extends to an entire function of complex spacetime, and translation by \(-i\mathbf q\) commutes with \(\mathcal P\). The first fact gives \(D_{\ell\rho}\widetilde f=0\), and together with the second it gives \(\operatorname{div}L_\ell=0\) in Section~\ref{sec:conservation}. Neither fact is available in this form for a variable-coefficient equation. The Gaussian filter also relies on the exact square identity for the matrix symbol in Section~\ref{sec:symbol}, and the \(tt\) entry \eqref{eq:time-component} is an exact Fourier multiplier applied to the wave. The recovery step in Section~\ref{sec:one-fiber} inverts that multiplier.

\subsection{The quaternionic Hessian and its cofactor}\label{sec:hessian}

The preceding construction gives a positive real tensor with divergence-free
columns. We now need a Hessian that pairs with this tensor and has a determinant
of the right size. The quaternionic Hessian serves both purposes. Its pairing
with a current is an ordinary real divergence-form operator, while its
determinant has one factor per quaternionic coordinate.

Why not the real or the complex Monge--Ampère operator? The current decides. Over each physical coordinate it has three real slots: the coordinate itself and the two real components of an auxiliary coordinate in \(\mathbb C_k\). (The auxiliary coordinates are complex because the two real states were paired as \(w_1+kw_2\).) Writing \(z_a=x_a+P_ai\) with \(P_a\in\mathbb C_k\), as in Lemma~\ref{lem:hessian} below, these slots are spanned by \(1\), \(i\) and \(j=ki\). Three real slots carry no complex structure, but they sit inside the quaternions, and the real tensor of Proposition~\ref{prop:wave-current} is the real representation of a quaternionic Hermitian matrix. The determinant adapted to this structure, the Moore determinant below, has one factor per quaternionic coordinate, that is, one per physical coordinate rather than one per real variable. On a fiber it therefore has degree \(n\), and this count is what produces \(p_c\) in the computation sketched in Section~\ref{sec:annular-idea}.

\subsubsection{The encoded gradient and Hessian.}
We first fix the conventions for linear algebra over the quaternions. Write
\((e_0,e_1,e_2,e_3)=(1,i,j,k)\). Conjugation reverses products,
\(\overline{pq}=\overline q\,\overline p\), and
\(\overline q q=|q|^2\). The real part satisfies
\(\operatorname{Re}(pq)=\operatorname{Re}(qp)\), as is seen by multiplying
out the real and imaginary parts. We treat quaternionic columns as a
right vector space: scalar multiplication takes \(v\) to
\(vq\), and a matrix \(A\) acts on the left, with
\((Av)_a=\sum_b A_{ab}v_b\), so that \(A(vq)=(Av)q\) by associativity.
The inner product of two columns is \(v^*w=\sum_a\overline{v_a}w_a\), and
\(\operatorname{Re}(v^*w)\) is the Euclidean inner product of their real
coordinates. These conventions fix the order of the factors in an outer
product \(vv^*\): its \((a,b)\) entry is \(v_a\overline{v_b}\), and
\(w^*vv^*w=|v^*w|^2\). In particular \(vv^*\) is positive semidefinite and,
for \(v\ne0\), has quaternionic rank one.
A Hermitian matrix \(A\) defines the form \(v^*Aw\), and its quadratic
form \(v^*Av\) is real because it equals its own conjugate. We write
\(A\geq0\) when this form is nonnegative and \(A>0\) when it is strictly
positive on every nonzero \(v\). A quaternionic unitary matrix \(U\) satisfies
\(U^*U=I\), which means that \(Uv\) and \(Uw\) have the same inner product
as \(v\) and \(w\).

The real representation~\eqref{eq:real-representation} expresses the
quadratic form of \(A\) in real coordinates. If
\(\widehat\xi_a=\sum_\alpha e_\alpha\xi_{a\alpha}\) encodes a real vector
\(\xi\) in three slots, expansion gives
\[
  \widehat\xi^*A\widehat\xi=\xi^{\mathsf t}A^{\mathbb R}\xi .
\]
With all four slots, the formula~\eqref{eq:real-representation} represents
\(A\) as a real linear operator on the whole quaternionic vector space.

The wave current has three real slots under each coordinate index, and
we package them as \(z_a=z_{a0}+iz_{a1}+jz_{a2}\). The proof of the
Dirichlet theorem will also use a fourth slot \(kz_{a3}\).
In the following formulas the indices \(\alpha,\beta\) range over the
slots present. For a real function \(u\), its encoded gradient and Hessian
are defined by
\[
  G_a(u)=\sum_\alpha e_\alpha\partial_{a\alpha}u,
  \qquad
  H_{ab}(u)=\sum_{\alpha,\beta}
      e_\alpha\overline{e_\beta}\,\partial_{a\alpha}\partial_{b\beta}u.
\]
Thus \(G(u)\) is the ordinary gradient with its entries grouped by coordinate
index, and \(H(u)\) is a linear function of the ordinary real Hessian. For a single coordinate with three slots \((x,p_1,p_2)\) it is the scalar \(u_{xx}+u_{p_1p_1}+u_{p_2p_2}\). The condition \(H(u)\geq0\) is therefore weaker than convexity: for example, \(u=-x^2+p_1^2+p_2^2\) has \(H(u)=2\), although its real Hessian is indefinite.

\begin{lemma}\label{lem:hessian}
The Hessian \(H(u)\) is Hermitian and \(\operatorname{tr}H(u)=\Delta u\).
It takes a real rank-one Hessian \(\xi\xi^{\mathsf t}\) to
\(\widehat\xi\widehat\xi^*\), where
\(\widehat\xi_a=\sum_\alpha e_\alpha\xi_{a\alpha}\) encodes the real covector
\(\xi\). In particular it preserves positivity.
It satisfies the product rule
\[
  H(uv)=uH(v)+vH(u)+G(u)G(v)^*+G(v)G(u)^*.
\]
It is covariant under quaternionic unitary changes of coordinates in four
slots and under real orthogonal changes of the coordinate indices in three slots.
If a function depends only on the complex coordinates \(P_a\in\mathbb C_k\)
in \(z_a=x_a+P_ai\), then
\[
    G_a(u)=2(\partial_{\overline P_a}u)i,
    \qquad H_{ab}(u)=4\partial_{\overline P_a}\partial_{P_b}u.
\]
In particular, if \(\Psi\) is a holomorphic map of the \(P\) variables with
complex derivative \(D\Psi\), then \(H(|\Psi|^2)=4(D\Psi)^*D\Psi\).
For any smooth real function of \(\Psi\), its encoded gradient lies in the
quaternionic image of \((D\Psi)^*\).
\end{lemma}

\begin{proof}
Taking the conjugate of \(H_{ab}(u)\) reverses its two quaternionic factors.
Interchanging the slot indices then gives \(H_{ba}(u)\). For \(a=b\), the
terms with unequal slots cancel in pairs by
\(e_\alpha\overline{e_\beta}+e_\beta\overline{e_\alpha}=0\) for
\(\alpha\ne\beta\). The terms with equal slots give
\(H_{aa}(u)=\sum_\alpha\partial_{a\alpha}^2u\), so summing over \(a\) gives
the real Laplacian.

If the real Hessian is \(\xi\xi^{\mathsf t}\), the \((a,b)\) entry factors as
\[
 \sum_{\alpha,\beta}
    e_\alpha\overline{e_\beta}\,\xi_{a\alpha}\xi_{b\beta}
 =\left(\sum_\alpha e_\alpha\xi_{a\alpha}\right)
  \left(\sum_\beta\overline{e_\beta}\xi_{b\beta}\right).
\]
The real coefficients commute with every quaternion, so this is exactly
\(\widehat\xi_a\widehat\xi_b^*\). Each such outer product is positive
because \(v^*\widehat\xi\widehat\xi^*v=|\widehat\xi^*v|^2\).
Decomposing a positive real Hessian into positive rank-one forms proves
the positivity assertion. The two mixed terms from the real product rule
factor in the same order to give the displayed formula for \(H(uv)\).

For covariance in four slots, make the change \(z=Uz'\) with \(U^*U=I\).
The linear function \(\operatorname{Re}(\widehat\xi^*z)\) now has encoded
covector \(U^*\widehat\xi\). Thus the rank-one formula transforms by
conjugation with \(U^*,U\). Real rank-one forms span the real symmetric
matrices, so this proves
\[
 H(u\circ U)(z')=U^*H(u)(Uz')U.
\]
The argument in three slots is the same when \(U\) is real orthogonal on
the coordinate index. Such a \(U\) preserves the three-slot subspace.

For the complex formulas, write \(P_a=p_a+kq_a\). Then
\(P_ai=ip_a+jq_a\), and the complex derivatives use
\(2\partial_{\overline P_a}=\partial_{p_a}+k\partial_{q_a}\) and
\(2\partial_{P_a}=\partial_{p_a}-k\partial_{q_a}\).
The gradient identity follows from \(ki=j\). For the Hessian, the four
products of the two remaining basis elements with their conjugates are
\(i\overline i=j\overline j=1\), \(i\overline j=-k\), and
\(j\overline i=k\). Consequently
\[
 H_{ab}(u)=
   \partial_{p_a}\partial_{p_b}u-k\partial_{p_a}\partial_{q_b}u
   +k\partial_{q_a}\partial_{p_b}u+\partial_{q_a}\partial_{q_b}u.
\]
This is the stated complex Hessian. For example, for \(u=|P_a|^2\)
the only nonzero entry is \(H_{aa}(u)=4\), the sum of the two real second
derivatives in \(p_a,q_a\).

For a holomorphic \(\Psi\), applying \(\partial_{P_b}\) to
\(\sum_\ell\overline{\Psi_\ell}\Psi_\ell\) differentiates only \(\Psi_\ell\).
Applying \(\partial_{\overline P_a}\) then differentiates only
\(\overline{\Psi_\ell}\). Hence
\[
 H_{ab}(|\Psi|^2)
   =4\sum_\ell \overline{\partial_{P_a}\Psi_\ell}\,\partial_{P_b}\Psi_\ell.
\]
For a real function \(\rho\) on the target of \(\Psi\), the same complex
chain rule gives
\(
 G(\rho\circ\Psi)=
 2(D\Psi)^*((\partial_{\overline\Psi_\ell}\rho)\circ\Psi)_\ell\, i.
\)
The order of the final factor \(i\) comes from the gradient formula.
\end{proof}

\subsubsection{The determinant and the normal cofactor.}
For Hermitian matrices write \(A:B=\operatorname{Re}\operatorname{tr}(AB)\).
The real part makes this pairing symmetric. A Hermitian quaternionic matrix
has real eigenvalues and a quaternionic orthonormal eigenbasis; the proof
of Lemma~\ref{lem:cofactor} recalls why. Its \emph{Moore determinant} is
the product of these eigenvalues with
multiplicity~\cite[Section~1.1]{alesker-linear}. For instance, for
\(c_1,c_2\in\mathbb R\) and \(q\in\mathbb H\),
\[
 \det\begin{pmatrix}c_1&q\\\overline q&c_2\end{pmatrix}=c_1c_2-|q|^2 .
\]
Throughout, the determinant of a Hermitian quaternionic matrix means its
Moore determinant. It is related to real determinants by a count of
eigenvalues: a quaternionic eigenvector with real eigenvalue \(\lambda\)
contributes a single factor \(\lambda\) to the Moore determinant, while
its quaternionic line has real dimension four and contributes four
factors of \(\lambda\) to the determinant of the four-slot real
representation. The determinant of the four-slot real representation is
therefore the fourth power of the Moore determinant.

For the Gram matrix \((D\Psi)^*D\Psi\) above, all entries lie in
\(\mathbb C_k\). Its complex unitary diagonalization is also a
quaternionic unitary diagonalization. Thus if \(D\Psi\) is square of
size \(m\), the determinant computed over either field gives
\(\det H(|\Psi|^2)=4^m|\det_{\mathbb C_k}D\Psi|^2\).
The fiber domains of Section~\ref{sec:fiber-potentials} are defined by a
squared norm of this type.

The cofactor \(C(A)\) describes the first-order change in the determinant
when \(A\) is perturbed. In particular \(v^*C(A)v\) measures the change
when one adds \(tvv^*\). In an orthonormal basis with first column \(v\),
this changes only the first diagonal entry. For a unit column \(v\), the
\emph{compression} \(A|_{v^\perp}\) is the matrix of the Hermitian form of
\(A\) on the quaternionic orthogonal complement \(v^\perp=\{w:v^*w=0\}\),
obtained by deleting the first row and column in such a basis. The
complement need not be invariant under \(A\). For the \(2\times2\) matrix
displayed above, the definition in Lemma~\ref{lem:cofactor} gives the
familiar formula
\begin{equation}\label{eq:cofactor-2x2}
 C\begin{pmatrix}c_1&q\\\overline q&c_2\end{pmatrix}
   =\begin{pmatrix}c_2&-q\\-\overline q&c_1\end{pmatrix}.
\end{equation}

\begin{lemma}\label{lem:cofactor}
On \(m\times m\) Hermitian quaternionic matrices the Moore determinant is
a homogeneous polynomial of degree \(m\), with cofactor \(C(A)\) determined by
\[
  C(A):B=\left.\frac{d}{dt}\det(A+tB)\right|_{t=0}.
\]
For \(A>0\), \(C(A)=(\det A)A^{-1}\) and \(\log\det A\) is concave and
increasing in the matrix order. If \(A\ge 0\), then \(C(A)\ge 0\), and if
\(A\) is also singular, \(C(A)\) kills its image.

For any unit column \(v\), and any Hermitian \(A\) and real \(t\),
\[
  v^*C(A)v=\det(A|_{v^\perp}),\qquad
  \det(A+tvv^*)=\det A+t\,v^*C(A)v.
\]
If
\(M=\begin{pmatrix}b&c^*\\c&D\end{pmatrix}\) with \(D>0\), and \(e\) is the
first unit column, then
\[
  \det M=(\det D)(b-c^*D^{-1}c),
  \qquad M-(b-c^*D^{-1}c)ee^*\ge 0.
\]
\end{lemma}

\begin{proof}
First, the usual Rayleigh quotient proof of the spectral theorem works with
the real inner product \(\operatorname{Re}(v^*w)\). Maximizing \(v^*Av\)
on the real unit sphere gives \(Av=\lambda v\) by first variation, with
\(\lambda=v^*Av\in\mathbb R\). The whole quaternionic line spanned by \(v\)
is an eigenspace because \(\lambda\) is real. Its quaternionic orthogonal
complement is invariant under \(A\), since \(v^*Aw=(Av)^*w=0\) there. Repeating
this argument gives a quaternionic orthonormal eigenbasis.
In that basis \(\operatorname{Re}\operatorname{tr}(A^\ell)\) is the sum of
the \(\ell\)-th powers of the eigenvalues. This equality is invariant under
change of basis, by cyclicity of the real part of the trace. Newton's
identities express their product as a polynomial in these power sums.
The determinant is therefore a polynomial in the real coefficients of \(A\).

For \(A>0\), the logarithm of the Moore determinant equals one quarter of
the real log determinant in four slots. Differentiating the real log
determinant along \(A+tB\) and returning to quaternionic notation gives
\[
 \begin{aligned}
 \left.\frac{d}{dt}\log\det(A+tB)\right|_{t=0}&=A^{-1}:B,\\
 \left.\frac{d^2}{dt^2}\log\det(A+tB)\right|_{t=0}
     &=-\operatorname{Re}\operatorname{tr}(A^{-1}BA^{-1}B).
 \end{aligned}
\]
Here we use that the real trace of a quaternionic linear matrix is four times the
real part of its quaternionic trace. For the first derivative, if \(B\geq0\)
then \(\operatorname{Re}\operatorname{tr}(A^{-1/2}BA^{-1/2})\geq0\).
For the second derivative,
\(\operatorname{Re}\operatorname{tr}(A^{-1}BA^{-1}B)\) is the sum of squares
of the eigenvalues of \(A^{-1/2}BA^{-1/2}\).
This proves monotonicity and concavity. Differentiating the determinant
itself gives \(C(A)=(\det A)A^{-1}\).

In an eigenbasis of \(A\), the cofactor thus replaces each diagonal
eigenvalue by the product of all the other eigenvalues. For \(A\geq0\),
including singular \(A\), this follows by letting the positive matrices
\(A+\epsilon I\) converge to \(A\). The cofactor is continuous
because it is the derivative of a polynomial. If \(A\) is singular, the
cofactor eigenvalue corresponding to each positive eigenvalue of \(A\)
contains a zero factor. This is why the cofactor kills the image of \(A\).
In particular, if \(A\) has a one-dimensional quaternionic kernel,
then its cofactor is supported on that line and its eigenvalue there is
the product of the positive eigenvalues of \(A\).

For the block formulas let \(q\in\mathbb H\) be the first coordinate of a
column and \(w\) the remaining coordinates. Completing the square gives
\[
 \begin{pmatrix}q\\w\end{pmatrix}^{\!*}
  M\begin{pmatrix}q\\w\end{pmatrix}
  =(w+D^{-1}cq)^*D(w+D^{-1}cq)
     +(b-c^*D^{-1}c)|q|^2.
\]
In particular, after subtracting the indicated multiple of \(ee^*\), the
form is nonnegative. The real change of variables
\((q,w)\mapsto(q,w+D^{-1}cq)\) is triangular with determinant one. Thus for
positive \(M\) the real determinant in four slots is
\((b-c^*D^{-1}c)^4(\det D)^4\). Taking the positive fourth root gives the
Schur determinant formula. The resulting identity is polynomial in \(b\)
and extends to every \(b\) for each fixed \(D>0\) and \(c\).

Finally, in a basis with first column \(v\), adding \(tvv^*\) changes
only \(b\), so the formula makes the determinant affine in \(t\) with
slope \(\det(A|_{v^\perp})\) when \(A>0\). Its derivative at \(t=0\) is
\(C(A):vv^*=v^*C(A)v\). Both identities extend to all Hermitian \(A\)
by polynomiality. For \(m=1\), the determinant of the empty compression
is understood to be \(1\), and the formulas follow from the scalar determinant.
\end{proof}

\subsubsection{Divergence and the boundary term.}
The formula for the real representation \(A^{\mathbb R}\) extends to four
slots by allowing \(e_3=k\). Using cyclicity of the real part and expanding
the definitions gives
\[
 \begin{aligned}
 A:H(u)
   &=\sum_{a,b,\alpha,\beta}
     \operatorname{Re}(\overline{e_\alpha}A_{ab}e_\beta)
         \,\partial_{a\alpha}\partial_{b\beta}u
     =A^{\mathbb R}:D^2u,\\
 (\nabla v)^{\mathsf t} A^{\mathbb R}\nabla u
   &=\operatorname{Re}\big(G(v)^*AG(u)\big).
 \end{aligned}
\]
In the first line \(D^2u\) and the pairing on the right are real. If the
real columns of \(A^{\mathbb R}\) have zero divergence, then the real
product rule gives \(A:H(u)=\operatorname{div}(A^{\mathbb R}\nabla u)\).
The divergence of the real \((b,\beta)\) column is
\[
  \operatorname{Re}\!\left(
    \left[\sum_{a,\alpha}\overline{e_\alpha}\,
                  \partial_{a\alpha}A_{ab}\right]e_\beta\right).
\]
Vanishing of the quaternionic sum in brackets therefore suffices for all the real
columns. Multiplication by a real scalar \(\chi_0\) preserves real
column conservation if \(AG(\chi_0)=0\). The extra term from
differentiating \(\chi_0\) is the real pairing of \(\nabla\chi_0\) with
that column, and \(AG(\chi_0)=0\) makes all such pairings vanish.

Alesker~\cite[Lemma~4.6]{alesker-cofactor} put the linearized quaternionic
Monge--Ampère operator in divergence form for positive definite Hessians.
The polynomial argument of OpenAI~\cite[Section~3.3]{openai} gives the
three-slot identity needed here, including for singular and indefinite
Hessians. The rank calculation underlying the argument is much like a
symbol calculation for a multilinear constant-coefficient operator.

\begin{lemma}[Cofactor divergence]\label{lem:cofactor-divergence}
For a smooth real function \(u\), \(C(H(u))\) satisfies
\[
   \sum_{a,\alpha}\overline{e_\alpha}\,
      \partial_{a\alpha} C(H(u))_{ab}=0
   \qquad\text{for every } b.
\]
\end{lemma}

For \(m=2\) the lemma can be checked by hand. By \eqref{eq:cofactor-2x2},
\(C(H(u))_{11}=H_{22}(u)\) and \(C(H(u))_{21}=-H_{21}(u)\), so for \(b=1\)
the claim is
\[
  \sum_\alpha\overline{e_\alpha}\,\partial_{1\alpha}H_{22}(u)
  =\sum_\alpha\overline{e_\alpha}\,\partial_{2\alpha}H_{21}(u).
\]
By the definition of the Hessian, the right side is
\(\sum_{\alpha,\beta,\gamma}\overline{e_\alpha}e_\beta\overline{e_\gamma}\,
\partial_{2\alpha}\partial_{2\beta}\partial_{1\gamma}u\). The derivatives are
symmetric in \(\alpha,\beta\), and
\(\overline{e_\alpha}e_\beta+\overline{e_\beta}e_\alpha=2\delta_{\alpha\beta}\),
so the right side equals
\(\sum_\gamma\overline{e_\gamma}\,\partial_{1\gamma}
\sum_\alpha\partial_{2\alpha}^2u\). This is the left side, because
\(H_{22}(u)=\sum_\alpha\partial_{2\alpha}^2u\). The case \(b=2\) is the
same with the two coordinates exchanged. The proof below treats all \(m\)
at once by polarization.

\begin{proof}
For \(m=1\) the cofactor is constant. For \(m>1\), take the coefficient of
\(t_1\cdots t_{m-1}\) in \(C(\sum_{\ell=1}^{m-1}t_\ell A_\ell)\), and call it
\(\mathcal C(A_1,\ldots,A_{m-1})\). Since \(C\) is homogeneous of degree
\(m-1\), this coefficient is symmetric and multilinear in the matrices.
When all of them equal \(A\), it is \((m-1)!C(A)\).

Take \(A_\ell=v_\ell v_\ell^*\) and \(t_\ell>0\). The matrix \(\sum_\ell t_\ell A_\ell\) has
kernel \(\{w:v_\ell^*w=0\ \text{for all }\ell\}\), since its quadratic form is
\(\sum_\ell t_\ell|v_\ell^*w|^2\). It is positive semidefinite with image the
quaternionic span of the \(v_\ell\), hence has rank at most \(m-1\).
Lemma~\ref{lem:cofactor} says its cofactor kills every \(v_\ell\). The
identity
\(C(\sum_\ell t_\ell A_\ell)v_{\ell'}=0\) holds on the open set where all \(t_\ell>0\). Taking its
\(t_1\cdots t_{m-1}\) coefficient gives
\(\mathcal C(A_1,\ldots,A_{m-1})v_{\ell'}=0\) for each \(\ell'\).

Now test the polarized differential identity on
\(u_\ell(z)=e^{\xi_\ell\cdot z}\) with real covectors \(\xi_\ell\), where the
dot product uses the real coordinates of \(z\).
By Lemma~\ref{lem:hessian},
\(H(u_\ell)=u_\ell\widehat\xi_\ell\widehat\xi_\ell^*\).
Multilinearity of \(\mathcal C\) factors out \(\prod_\ell u_\ell\).
Differentiating this product and canceling it from the divergence, the
quaternionic divergence is then the row
\[
  \left(\sum_\ell\widehat\xi_\ell\right)^*
  \mathcal C(\widehat\xi_1\widehat\xi_1^*,\ldots,
              \widehat\xi_{m-1}\widehat\xi_{m-1}^*).
\]
The matrix is Hermitian and kills each \(\widehat\xi_\ell\), so the displayed
row vanishes. To see why these tests suffice, expand the polarized
differential expression in products of real derivatives
\(\prod_\ell\partial^{\gamma_\ell}u_\ell\). The coefficients are constant.
Evaluation on the exponentials replaces each
\(\partial^{\gamma_\ell}u_\ell/u_\ell\) by \(\xi_\ell^{\gamma_\ell}\).
The resulting polynomial is zero for all \(\xi_\ell\), so all its
coefficients vanish. This proves the polarized identity on smooth
functions. Setting all \(u_\ell=u\) and dividing by \((m-1)!\) proves the
cofactor identity.
\end{proof}

For a smooth domain \(\Omega=\{r<0\}\), write
\(\mathbf n=\nabla r/|\nabla r|\) for the outward real unit normal and
\(\nu=G(r)/|\nabla r|\) for its quaternionic encoding. The quaternionic
orthogonal complement \(\nu^\perp=\{w:\nu^*w=0\}\) is in general smaller
than the real tangent hyperplane, which imposes only real orthogonality
to \(\mathbf n\). If \(A\) has
divergence-free real columns, apply the real divergence theorem to
\(v A^{\mathbb R}\nabla u-u A^{\mathbb R}\nabla v\). The two terms
with a derivative on \(u\) and on \(v\) cancel by symmetry of
\(A^{\mathbb R}\). Encoding the two normal components gives
\begin{equation}\label{eq:green}
  \int_\Omega (vA:H(u)-uA:H(v))
  =\int_{\partial\Omega}
     \operatorname{Re}\big(v\nu^*AG(u)-u\nu^*AG(v)\big)\,dS .
\end{equation}

Suppose that \(u\) vanishes on the boundary. Taking \(r\) as the last local
coordinate and integrating the derivative of \(u\) in that coordinate
from \(0\) to \(r\) shows that \(u=\rho r\) for a smooth \(\rho\). On the
boundary \(G(u)=\rho G(r)=\nu\partial_{\mathbf n}u\), so the first normal
term in \eqref{eq:green} is \(v(\nu^*A\nu)\partial_{\mathbf n}u\).
The product rule gives on the boundary
\[
 H(u)=\rho H(r)+G(\rho)G(r)^*+G(r)G(\rho)^* .
\]
On vectors in \(\nu^\perp\) the last two terms vanish as a quadratic form,
because \(G(r)=|\nabla r|\nu\). Equivalently, the compression satisfies
\begin{equation}\label{eq:boundary-compression}
  H(u)|_{\nu^\perp}=\rho H(r)|_{\nu^\perp}
  \qquad\text{on }\partial\Omega .
\end{equation}
For \(A=C(H(u))\), the coefficient \(\nu^*A\nu\) is therefore the
determinant of this compressed Hessian, by Lemma~\ref{lem:cofactor}. The
energy argument will use the same calculation with a background potential
added to \(u\). A null direction of the compressed Hessian then makes the
corresponding boundary term in Green's formula vanish.

\subsection{The periodic Dirichlet problem}

On a fiber we will solve the quaternionic Monge--Ampère equation
\(\det H(B_0+u)=g\), where \(g\) is the prescribed density and \(B_0\) is a
smooth background potential. The correction \(u\) vanishes at the boundary,
so the total potential \(B_0+u\) takes the boundary values \(B_0\). We call
the total potential \emph{admissible} if its quaternionic Hessian is positive
definite. This positivity is what makes the determinant equation elliptic.
(A smooth function of four-slot variables is called quaternionic
plurisubharmonic when its quaternionic Hessian is
nonnegative~\cite{alesker-linear}; in the three-slot setting one may test
this on its constant extension in the missing slot. Admissibility requires
the strict inequality.)
For one quaternionic coordinate, \(H(\psi)=\Delta\psi\), so the equation
is a Poisson equation. In higher dimensions the determinant fixes only
the product of the eigenvalues, and the existence proof must control
each eigenvalue. We will use two consequences of the theory: the maximum
principle lets us compare the depth \(-u\) with explicit functions, and
smoothness up to the moving boundary lets us differentiate the energy
in the three base variables.

Alesker~\cite[Theorem~1.4]{alesker-dirichlet} proved smooth Dirichlet
solvability for this equation on quaternionic balls. We give the
periodic, three-slot version of OpenAI~\cite[Appendix~A]{openai}.

\begin{theorem}[Periodic Dirichlet problem]\label{thm:dirichlet}
Let \(\Lambda\subset\mathbb R^m\) be a full lattice, let \(V\subset\mathbb R^{2m}\)
be a smooth bounded domain, and put \(\Omega=(\mathbb R^m/\Lambda)\times V\).
Use the periodic coordinates in the real slots of the Hessian \(H\), and the
coordinates of \(V\) in its two imaginary slots. Suppose that \(r\in C^\infty(\overline V)\)
is a defining function for \(V\), meaning that \(r<0\) in \(V\), \(r=0\) on
\(\partial V\) and \(dr\ne0\) on \(\partial V\). Assume \(H(r)>0\) on
\(\overline\Omega\). Let \(B_0,g\in C^\infty(\overline\Omega)\) be real and satisfy
\(H(B_0)\geq0\) and \(g>\det H(B_0)\) on \(\overline\Omega\). Then there is a unique
\(u\in C^\infty(\overline\Omega)\) satisfying
\[
  \det H(B_0+u)=g,\qquad H(B_0+u)>0,\qquad u|_{\partial\Omega}=0.
\]
For some \(b\geq b'>0\) one has \(br\leq u\leq b'r\). For smooth
finite-dimensional families of such data and domains the constants can be
chosen locally uniformly, and the solutions depend smoothly on the
parameters up to the boundary in any smooth local identification with a
fixed domain.
\end{theorem}

The problem is periodic in the torus factor only; the boundary of \(\Omega\) is \((\mathbb R^m/\Lambda)\times\partial V\). The functions \(br\) and \(b'r\) are lower and upper barriers for \(u\): since \(r<0\) inside, the larger constant gives the lower function.

The later sections use Theorem~\ref{thm:dirichlet} only through the following properties: existence and uniqueness of the admissible solution; the barrier bounds \(br\leq u\leq b'r\), in particular \(u\leq0\); smooth dependence on parameters up to the boundary, which lets us differentiate the energy in Lemma~\ref{lem:energy-balance}; and the comparison principle below, on which the two comparison arguments in Proposition~\ref{prop:fiber-estimates} rest. The proof of Theorem~\ref{thm:dirichlet}, given in Appendix~\ref{app:dirichlet}, adapts the method of Caffarelli, Nirenberg and Spruck~\cite{caffarelli-nirenberg-spruck}.

\begin{lemma}[Comparison]\label{lem:comparison}
Let \(\psi_1,\psi_2\) be smooth potentials that are admissible on \(\overline\Omega\). If \(\det H(\psi_1)\geq\det H(\psi_2)\) in \(\Omega\) and \(\psi_1\leq\psi_2\) on \(\partial\Omega\), then \(\psi_1\leq\psi_2\) in \(\Omega\). The same holds with \(\Omega\) replaced by an open subset, such as a ball inside a fiber.
\end{lemma}

\begin{proof}
For potentials \(\psi_1,\psi_2\) that are admissible on the closure,
the entire segment between their Hessians is positive. Applying the fundamental theorem of
calculus and the derivative of the log determinant from
Lemma~\ref{lem:cofactor} gives
\[
 \begin{aligned}
 &\log\det H(\psi_1)-\log\det H(\psi_2)\\
 &\quad=
 \left(\int_0^1 \big((1-\theta)H(\psi_2)+\theta H(\psi_1)\big)^{-1}
                   \,d\theta\right):H(\psi_1-\psi_2).
 \end{aligned}
\]
For a Hermitian coefficient \(A>0\), the real operator \(w\mapsto A:H(w)\)
has symbol \(\widehat\xi^* A\widehat\xi\) by the rank-one identity in
Lemma~\ref{lem:hessian}. Encoding the covector \(\xi\) in its slots
preserves its Euclidean length. Thus eigenvalue bounds \(cI\leq A\leq CI\)
give the bounds \(c|\xi|^2\leq\widehat\xi^* A\widehat\xi
\leq C|\xi|^2\) needed for uniform ellipticity. The averaged inverse
in the display has such bounds by compactness. If
\(\det H(\psi_1)\geq\det H(\psi_2)\), then \(\psi_1-\psi_2\) is subharmonic
for this real elliptic operator. Hence \(\psi_1\leq\psi_2\) on the
boundary implies \(\psi_1\leq\psi_2\) throughout. This gives uniqueness
by applying comparison in both directions. The argument also works
on subdomains such as the comparison balls inside a fiber.
\end{proof}

\section{Proof of the local smoothing conjecture}\label{sec:proof}

\subsection{The fibers and their Dirichlet potentials}\label{sec:fiber-potentials}

Work on the fixed spatial torus with annular data in the cap \(|v|\leq\gamma\), where \(\gamma>0\) will be chosen small. We may assume \(\|f\|_2=1\). Set \(a=K/R\), where \(K\) is a large fixed constant to be chosen after \(\gamma\), and assume \(4a<1\). The base point \(z=(t,\tau)\) has three real coordinates, with \(\tau\in\mathbb C_k\) and \(a<s=|\tau|<4\). The \emph{fiber} over \(z\) is the set of points \((x,P,Z)\) above \(z\): the spatial point \(x\) runs over the torus, and \((P,Z)\) are the \(2n\) real auxiliary coordinates attached to the spatial coordinates. In the fiber we consider the domain
\[
 \Omega_z=\{r<0\},\qquad
 r=100|W|^2+100|Q|^2-1,\qquad Q=\Phi/a ,
\]
in which \(x\) is unrestricted. Here \(r\) is a defining function: its negative set is the domain, and its zero set, on which the gradient is nonzero, is the boundary of \(\Omega_z\). We call this boundary the \emph{fiber wall}. The full quaternionic Hessian has size \((n+1)\times(n+1)\), while its spatial block at fixed \(z\) has size \(n\times n\). We write \(H_{\mathrm{fib}}\) and \(G_{\mathrm{fib}}\) for the Hessian and gradient at fixed \(z\). Only these two objects omit the base derivatives. The current is conserved in all variables, and on a fiber the spatial point \(x\) is an integration variable like \(P\) and \(Z\).

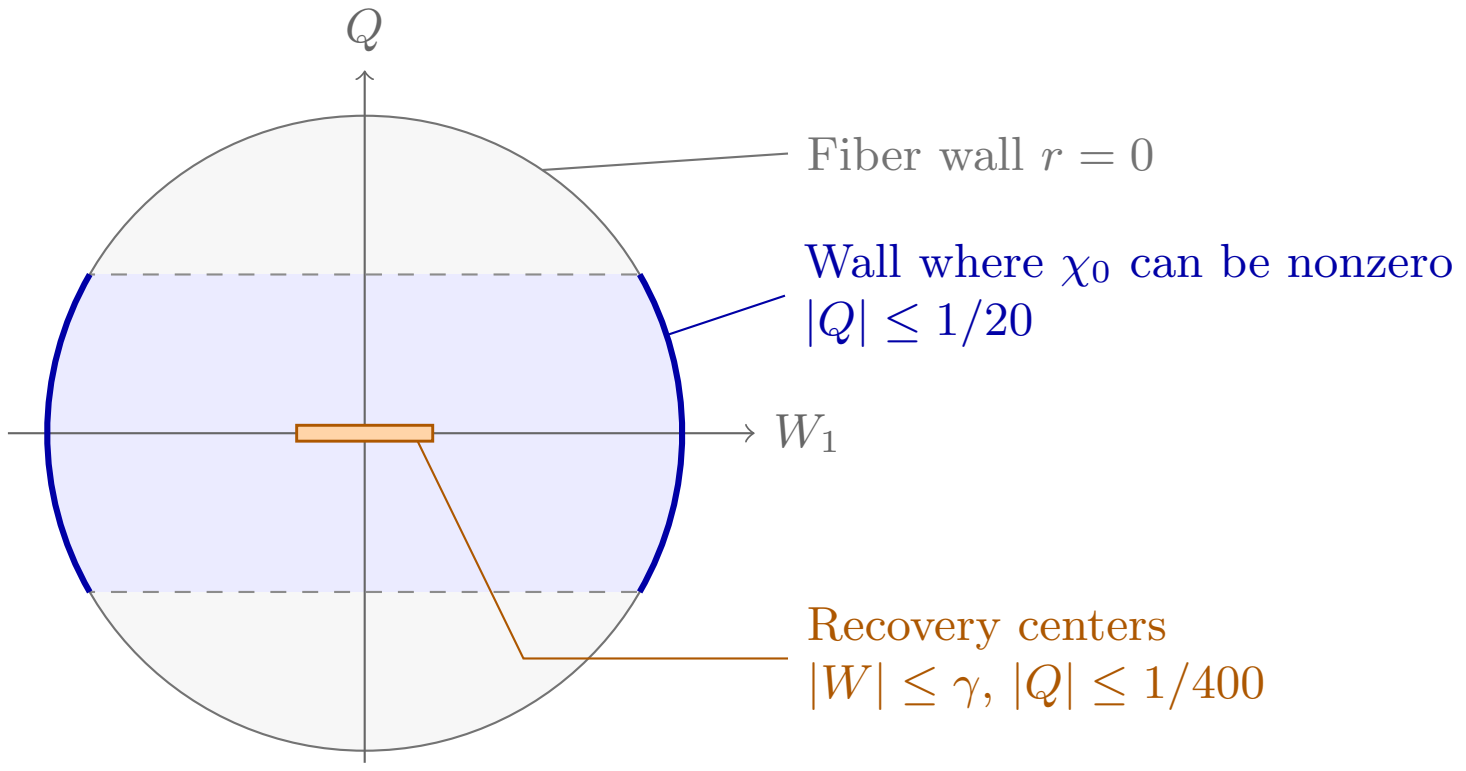
\begin{figure}[tbp]
\centering
\begin{tikzpicture}[font=\small]
  \fill[black!3] (0,0) circle[radius=2.1];
  \begin{scope}
    \clip (0,0) circle[radius=2.1];
    \fill[blue!8] (-2.1,-1.05) rectangle (2.1,1.05);
    \draw[dashed,black!45] (-2.1,1.05) -- (2.1,1.05);
    \draw[dashed,black!45] (-2.1,-1.05) -- (2.1,-1.05);
  \end{scope}
  \draw[black!55] (0,0) circle[radius=2.1];
  \draw[->,black!60] (-2.36,0) -- (2.58,0) node[right] {\(W_1\)};
  \draw[->,black!60] (0,-2.18) -- (0,2.40) node[above] {\(Q\)};
  \draw[blue!65!black,line width=1.15pt]
    (-30:2.1) arc[start angle=-30,end angle=30,radius=2.1];
  \draw[blue!65!black,line width=1.15pt]
    (150:2.1) arc[start angle=150,end angle=210,radius=2.1];
  \filldraw[fill=orange!34,draw=orange!68!black,line width=.55pt]
    (-.45,-.0525) rectangle (.45,.0525);
  \draw[black!55] (56:2.1) -- (2.80,1.85)
    node[right,align=left] {Fiber wall \(r=0\)};
  \draw[blue!65!black] (18:2.1) -- (2.78,.91)
    node[right,align=left] {Wall where \(\chi_0\)\ can be nonzero\\
                             \(|Q|\leq 1/20\)};
  \draw[orange!68!black] (.35,-.0525) -- (1.05,-1.49) -- (2.80,-1.49)
    node[right,align=left] {Recovery centers\\
                             \(|W|\leq\gamma\), \(|Q|\leq 1/400\)};
\end{tikzpicture}
\caption{A real slice of the fiber at fixed \((t,\tau,x)\), with \(W_1\) and \(Q\)
real and the remaining entries of \(W\) zero. The blue arcs have
\(|W|\geq\sqrt{3}/20\), so the Gaussian is small there for the chosen cap.
The orange strip marks the centers used to recover the wave. Its
horizontal width is schematic.}
\label{fig:fiber}
\end{figure}

The coordinates \(W\) and \(Q\) have separate roles. The variable \(W\) is the center of the Gaussian in the cap variable \(v\). The coordinate \(\Phi\) is constant along the directions of the current, so we can cut off in \(Q=\Phi/a\) while retaining the conservation law. The cutoff confines the portion of the wall carrying wave flux to \(|W|\geq\sqrt{3}/20\), separated from \(|v|\leq\gamma\). The recovery region has \(|W|\leq\gamma\), where the Gaussian is large enough to invert at \(s=2a\). Figure~\ref{fig:fiber} shows the two regions.

This separation is easy to see in the coordinates \((W,Q)\). Since
\(Z=\sqrt{2}Y-\tau\) and \(P=WY\), the equation for \(\Phi\) reads
\(aQ=Y(1+W\cdot W/2)-\sqrt{2}\,\tau\), so
\[
 Y=\frac{\sqrt{2}\,\tau+aQ}{1+W\cdot W/2},
 \qquad P=WY,\qquad Z=\sqrt{2}\,Y-\tau .
\]
This chart identifies the \((P,Z)\) part of every domain \(\Omega_z\) with the same complex ball in \((W,Q)\). The denominator is bounded away from zero since \(|W|\leq 1/10\), and \(|aQ|\leq a/10<s/10\), so \(|Y|\asymp s\). Integration on \(\Omega_z\) uses \(dx\,dP\,dZ\) in the original real coordinates, with two real dimensions for each \(\mathbb C_k\) coordinate. The determinant is likewise taken in the original coordinates \((x,P,Z)\) and not in the chart \((W,Q)\).

\begin{lemma}\label{lem:fiber-geometry}
The spatial Hessian \(H_{\mathrm{fib}}(r)\) is positive definite, and
\[
 \begin{gathered}
 H_{\mathrm{fib}}(r)\gtrsim s^{-2}I,\qquad
 G_{\mathrm{fib}}(r)G_{\mathrm{fib}}(r)^*\leq H_{\mathrm{fib}}(r),\\
 dP\,dZ\asymp a^2s^{2d}\,dW\,dQ,\qquad
 \det H_{\mathrm{fib}}(r)\asymp a^{-2}s^{-2d}.
 \end{gathered}
\]
The full Hessian \(H(r)\) has a one-dimensional quaternionic kernel contained in \(G(r)^\perp\) and \(\ker H(|\Phi|^2)\). Moreover, \(H_{\mathrm{fib}}(|\Phi|^2)\) has rank one with its positive eigenvalue comparable to one, and \(|G(|W|^2)|\lesssim s^{-1}\).
\end{lemma}

\begin{proof}
The full coordinate change from \((\tau,P,Z)\) to \((\tau,W,Q)\)
is invertible on the region under consideration. Discarding the first
output leaves the complex derivative \(D(10W,10Q)\) of rank \(n\).
Lemma~\ref{lem:hessian} gives
\[
 H(r)=4D(10W,10Q)^*D(10W,10Q),\qquad
 G(r)=2D(10W,10Q)^*(10W,10Q)^{\mathsf t}i.
\]
The complex rank is also the quaternionic rank, by row reduction over
\(\mathbb C_k\). The kernel is a quaternionic line, represented by
varying \(\tau\) with \((W,Q)\) fixed. For a vector on that line,
both its inner product with \(G(r)\) and its image under the last row
of \(D(10W,10Q)\) vanish. Since \(\Phi=aQ\), the Hessian formula
for \(|\Phi|^2\) puts the same line in \(\ker H(|\Phi|^2)\).

On a fixed fiber the inverse chart gives
\[
 dY=\frac{a\,dQ-YW\cdot dW}{1+W\cdot W/2},\qquad
 dP=Y\,dW+W\,dY,\qquad dZ=\sqrt{2}\,dY.
\]
The columns of its derivative in \(W\) have size \(O(s)\)
and the column in \(Q\) has size \(O(a)\). Thus the inverse
of \(D_{\mathrm{fib}}(10W,10Q)\) has norm \(O(s)\), since \(a<s\).
The formula for \(H(r)\) restricted to the fiber gives the desired
lower bound. For the gradient inequality, test against any quaternionic
column \(w\). Cauchy--Schwarz and \(|(10W,10Q)|\leq1\) give
\[
 |G_{\mathrm{fib}}(r)^*w|^2
 \leq 4|D_{\mathrm{fib}}(10W,10Q)w|^2
 =w^*H_{\mathrm{fib}}(r)w .
\]

To compute the volume, wedge \(dP_1,\ldots,dP_d\) with \(dZ\).
All terms involving \(dY\) from a factor \(dP_\ell\) vanish, since \(dZ\)
is a multiple of \(dY\). In the remaining expression only the \(dQ\)
part of \(dY\) survives. Hence
\[
 dP_1\wedge\cdots\wedge dP_d\wedge dZ
 =
 \frac{\sqrt{2}\,aY^d}{1+W\cdot W/2}\,
 dW_1\wedge\cdots\wedge dW_d\wedge dQ.
\]
The absolute square of the complex determinant is the real volume
Jacobian, so \(dP\,dZ\asymp a^2s^{2d}\,dW\,dQ\).
Also
\[
 \det H_{\mathrm{fib}}(r)
   =4^n|\det\nolimits_{\mathbb C_k}D_{\mathrm{fib}}(10W,10Q)|^2
   \asymp a^{-2}s^{-2d}.
\]

Finally, \(\partial_P\Phi=W\),
\(\partial_Z\Phi=(1-W\cdot W/2)/\sqrt{2}\), and
\(\partial_\tau\Phi=-(1+W\cdot W/2)/\sqrt{2}\). The fiber derivative
of \(\Phi\) is bounded above and away from zero. Thus
\(H_{\mathrm{fib}}(|\Phi|^2)=4(D_{\mathrm{fib}}\Phi)^*D_{\mathrm{fib}}\Phi\)
has one positive eigenvalue of size comparable to one.
Differentiating \(W=P/Y\) gives \(|DW|\lesssim s^{-1}\).
The gradient formula for \(|W|^2\) and the bound on \(|W|\) then give
\(|G(|W|^2)|\lesssim s^{-1}\).
\end{proof}

For \(\delta>0\), let \(\mathcal K_\delta\) be the periodization of the probability kernel
\[
 c_n\delta^{-n}(1+|x/\delta|^2)^{-n-2}.
\]
Take a smooth function \(0\leq\chi_0\leq1\) equal to one for \(|Q|\leq1/100\) and zero for \(|Q|\geq1/20\). Smooth the current of Proposition~\ref{prop:wave-current} in \(x\), cut it off in \(Q\), and add a small cofactor current:
\[
 A_w=\chi_0(Q)(\mathcal K_{1/R}*_x A_0),\qquad
 A=A_w+R^{-50n}C(H(r)),\qquad \mu=A_{tt}.
\]
We call \(\mu\) the \emph{density}. Convolution in \(x\) commutes with each derivative in the column
divergence. Also, \(DQ\) is independent of \(x\) and
\((DQ)A_0=0\) by Proposition~\ref{prop:wave-current}. Hence the
smoothed matrix annihilates the image of \((DQ)^*\), which contains
\(G(\chi_0(Q))\). The derivative of the cutoff in the product rule
for column divergence therefore gives zero.

The cofactor term supplies a strictly positive \(tt\) entry even at
points where the wave current vanishes. By
Lemmas~\ref{lem:cofactor} and~\ref{lem:cofactor-divergence},
\(C(H(r))\) is positive semidefinite with divergence-free real columns.
As a linear map it is supported on the one-dimensional kernel of
\(H(r)\), and its \(tt\) entry is the determinant of the spatial block,
\(\det H_{\mathrm{fib}}(r)>0\). Lemma~\ref{lem:fiber-geometry} puts
this kernel in both \(G(r)^\perp\) and \(\ker H(|\Phi|^2)\), so the
regularization contributes neither wall flux nor pairing with
\(H(|\Phi|^2)\). The wave term also has zero pairing with
\(H(|\Phi|^2)=4(D\Phi)^*D\Phi\), because \((D\Phi)A_w=0\).
In summary, \(A\) is positive semidefinite with divergence-free real
columns, \(\mu>0\), and \(A:H(|\Phi|^2)=0\).

Put \(B_0=\lambda|\Phi|^2\) with \(\lambda=R^{30n}\). The positive eigenvalue of \(H_{\mathrm{fib}}(B_0)\) is comparable to \(\lambda\). To bound \(\mu\), recall that Parseval and the count of lattice points in the annulus
give \(\sum_{\kappa,v}|c_{\kappa,v}|\lesssim R^{n/2}\). Each Gaussian
amplitude in the current is \(O(R^d)\), since \(s\geq a\geq 1/R\), so
\(\|A_0\|\lesssim R^{n+2d}\). The cofactor contributes at most
\(CR^{-50n}a^{-2}s^{-2d}\leq CR^{-50n+2+2d}\) to \(\mu\).
For sufficiently large \(R\), these two estimates give the density ceiling
\begin{equation}\label{eq:density-ceiling}
 \mu^{1/d}\leq\lambda .
\end{equation}

Apply Theorem~\ref{thm:dirichlet} on each fiber, with defining function \(r\), background \(B_0\) and prescribed density \(g=\lambda\mu\). The defining Hessian is strictly positive, the background Hessian \(H_{\mathrm{fib}}(B_0)\) is nonnegative with determinant zero, and the prescribed density is strictly positive. We obtain a solution, smooth up to the boundary and depending smoothly on \((t,\tau)\), of
\begin{equation}\label{eq:fiber-dirichlet}
 H_{\mathrm{fib}}(B_0+u)>0,\qquad
 \det H_{\mathrm{fib}}(B_0+u)=\lambda\mu,\qquad
 u|_{\partial\Omega_z}=0.
\end{equation}
The Dirichlet barriers give \(u\leq0\). We call \(-u\) the \emph{depth}: it
measures how far the potential lies below the background.

To see what depth to expect, diagonalize the rank-one matrix
\(H_{\mathrm{fib}}(B_0)\) and add \(\alpha I\). The \(d\) formerly null
eigenvalues become \(\alpha\), and the positive eigenvalue becomes
\(\beta+\alpha\), with \(\beta\asymp\lambda\). Thus
\begin{equation}\label{eq:background-determinant}
 \det\bigl(H_{\mathrm{fib}}(B_0)+\alpha I\bigr)
    =(\beta+\alpha)\alpha^d\asymp\lambda\alpha^d
    \qquad (0<\alpha\lesssim\lambda).
\end{equation}
A density of order \(\mu\) therefore asks for curvature of order
\(\mu^{1/d}\), which lies in this range by~\eqref{eq:density-ceiling}.
This is the role of the background. It takes up one of the \(n\) factors
of the determinant, at the fixed size \(\lambda\), and leaves \(d\) factors
for the curvature. Without it, the same count would give curvature of
order \(\mu^{1/n}\) and, in the energy, the power \(1+1/n\) of the density
in place of \(1+1/d=p_c/2\). The background is also invisible to the
current, since \(A:H(B_0)=0\). A
quadratic function with this curvature that vanishes on the boundary of a
ball of radius \(1/R\) has depth of order \(R^{-2}\mu^{1/d}\) at its
center. In Proposition~\ref{prop:fiber-estimates} we will use this
quadratic as a comparison function when we control the density from below
on the whole ball.

\subsection{The energy and its flux}\label{sec:energy}

The energy of a fiber pairs the depth with the density. The raw energy is
\[
 F_0(z)=\int_{\Omega_z}(-u)\mu\,dx\,dP\,dZ .
\]
By the discussion at the end of Section~\ref{sec:fiber-potentials} the depth
is governed by \(\mu^{1/d}\), so \(F_0\) measures the power \(1+1/d=p_c/2\)
of a squared wave. We want an inequality for the Laplacian of the energy in
the base variables.
Write \(\tau=\tau_1+k\tau_2\). The three real variables of
\(z=(t,\tau)\) carry the Laplacian
\[
 \Delta_z=\partial_t^2+\partial_{\tau_1}^2+\partial_{\tau_2}^2 .
\]
These variables are the three slots of the time coordinate, and the
subscript \(tt\) names the corresponding diagonal entry of a quaternionic
matrix. For example, if a test
function \(\eta\) depends only on \(z\), then \(H(\eta)\) has just
the entry \(H(\eta)_{tt}=\Delta_z\eta\).

The real column conservation of \(A\) lets us move these derivatives
from a test function to \(u\). As we show below, the interior term in
\(\Delta_z F_0\) is \(-\int_{\Omega_z} A:H(u)\,dx\,dP\,dZ\).
Pairing a positive divergence-free tensor with a Hessian has a precedent
in Serre's compensated integrability for real Monge--Ampère Hessians
\cite[Section~5.1]{serre-divergence}. Here the fiber Hessian is positive
definite. To extract the part of \(A:H(u)\) that has a sign, we use
a Schur complement as in Serre's diagonal-entry estimate
\cite[Proposition~2.2]{serre-mixed}.

Set \(D_\theta=H_{\mathrm{fib}}(B_0+\theta u)\) for \(0\leq\theta\leq1\).
The full Hessian has the block decomposition
\[
 \begin{gathered}
 H(B_0+u)=\begin{pmatrix}h&c^*\\c&D_1\end{pmatrix}
 = (h-c^*D_1^{-1}c)e_te_t^*+\Sigma,\qquad
 \Sigma=\begin{pmatrix}c^*D_1^{-1}c&c^*\\c&D_1\end{pmatrix}\geq0 .
 \end{gathered}
\]
Here \(e_t\) is the first coordinate column, and the real scalar \(h-c^*D_1^{-1}c\) is the \emph{Schur complement} of the spatial block \(D_1\). It is what remains of the \(tt\) entry after completing the square in the spatial variables, and it has no definite sign. The remainder \(\Sigma\), on the other hand, is positive semidefinite:
its quadratic form on a time component \(q\)
and a fiber column \(v\) is
\((v+D_1^{-1}cq)^*D_1(v+D_1^{-1}cq)\), which is nonnegative since
\(D_1>0\). Consequently \(A:\Sigma\geq0\), because in an orthonormal
eigenbasis of \(\Sigma\) the pairing
\(A:\Sigma\) is a sum of nonnegative eigenvalues of \(\Sigma\)
times nonnegative diagonal entries of \(A\).

By the Schur determinant formula in Lemma~\ref{lem:cofactor}, the scalar
\(h-c^*D_1^{-1}c\)
is \(\det H(B_0+u)/\det D_1\). Pairing \(e_te_t^*\) with
\(A\) gives \(A_{tt}=\mu\), and the equation \(\det D_1=\lambda\mu\)
then cancels \(\mu\). Since \(A:H(B_0)=0\), we obtain
\begin{equation}\label{eq:energy-schur}
 A:H(u)=A:\Sigma+\lambda^{-1}\det H(B_0+u).
\end{equation}
The first term on the right contributes with a favorable sign to
\(\Delta_z F_0\). The full determinant has no definite sign and
must be canceled.

To find the correction, consider the cofactor of \(H(B_0+\theta u)\)
for \(0\leq\theta\leq1\). Its \(tt\) entry is the determinant of the
fiber block \(D_\theta\), and its contraction against \(H(u)\) is
\(\partial_\theta\det H(B_0+\theta u)\).
If we start with a fiber integral of \(u\det D_\theta\), the
same Green identity that we use for \(F_0\) therefore produces a derivative of
the full determinant. Integrating in \(\theta\) gives the full
determinant at \(\theta=1\), since the rank-one background has zero
full determinant at \(\theta=0\). This motivates the correction
\[
 E(z)=\lambda^{-1}\int_{\Omega_z}
      u\int_0^1\det D_\theta\,d\theta\,dx\,dP\,dZ,
 \qquad F=F_0+E.
\]
Using \(\det D_1=\lambda\mu\), we can also write the corrected energy as
\begin{equation}\label{eq:energy-J}
 \lambda F=\int_{\Omega_z} (-u)
     \left(\det D_1-\int_0^1\det D_\theta\,d\theta\right)\,dx\,dP\,dZ .
\end{equation}

\begin{remark}[A model in two real variables]
The cancellation can be seen by hand in the simplest real analogue: one
base variable \(z\), one fiber variable \(x\in[-1,1]\), the ordinary
Hessian \(D^2\) in place of \(H\), no background, and \(\lambda=1\). Let
\(A\) be a smooth field of positive semidefinite symmetric \(2\times2\)
matrices with divergence-free columns and \(A_{zz}>0\), and let \(u(z,x)\)
solve the fiber equation \(u_{xx}=A_{zz}\) with \(u(z,\pm1)=0\). The fiber
Hessian of \(\theta u\) is the scalar \(\theta u_{xx}\), so the two parts of
the energy are
\[
 F_0(z)=\int_{-1}^1(-u)A_{zz}\,dx,\qquad
 E(z)=\int_{-1}^1u\int_0^1\theta u_{xx}\,d\theta\,dx
     =-\frac12\int_{-1}^1u_x^2\,dx .
\]
Differentiate \(E\) twice in \(z\) and integrate by parts in \(x\). The
boundary terms vanish because \(u_{zz}=0\) at \(x=\pm1\), and we find
\[
 E''(z)=\int_{-1}^1\bigl(u_{zz}u_{xx}-u_{zx}^2\bigr)\,dx
       =\int_{-1}^1\det D^2u\,dx .
\]
On the other hand, with \(e_z=(1,0)^{\mathsf t}\), the Schur decomposition is
\[
 D^2u=\frac{\det D^2u}{u_{xx}}\,e_ze_z^{\mathsf t}+\Sigma,
 \qquad
 \Sigma=\begin{pmatrix}u_{zx}^2/u_{xx}&u_{zx}\\u_{zx}&u_{xx}\end{pmatrix}\geq0,
\]
and pairing it with \(A\) gives \(A:D^2u=A:\Sigma+\det D^2u\), because
\(A_{zz}=u_{xx}\). Finally, since the columns of \(A\) are divergence free,
two integrations by parts give
\[
 F_0''(z)=-\int_{-1}^1A:D^2u\,dx+\bigl[A_{xx}u_x\bigr]_{x=-1}^{x=1}.
\]
Adding the three formulas, the determinant cancels:
\[
 (F_0+E)''+\int_{-1}^1A:\Sigma\,dx=\bigl[A_{xx}u_x\bigr]_{x=-1}^{x=1}.
\]
This is the model of the energy balance \eqref{eq:energy-balance} in
Lemma~\ref{lem:energy-balance} below, and its right side is the flux
through the wall \(x=\pm1\). The model has no room for the background or
for the exponent \(p_c\); it only shows the cancellation.
\end{remark}

\begin{remark}[Energies from Kähler geometry]
The segment integral in \(E\) is the quaternionic version of Mabuchi's energy
formula \cite[(2.2.1)]{mabuchi}. It is an antiderivative of the
determinant in the following sense. On a fixed fiber vary \(u\) in a
direction \(v\) that vanishes at the wall. The derivative of
\(\lambda E\) has one term from varying \(u\) outside the segment
integral and one term from varying \(D_\theta\). In the second term,
the cofactor has zero column divergence, so Green's formula moves
\(H_{\mathrm{fib}}\) from \(v\) onto \(u\). Both boundary values are
zero. The derivative becomes
\[
 \int_{\Omega_z} v\int_0^1
    \bigl[\det D_\theta+\theta\,C(D_\theta):H_{\mathrm{fib}}(u)\bigr]
       \,d\theta\,dx\,dP\,dZ
 =\int_{\Omega_z}v\det D_1\,dx\,dP\,dZ ,
\]
since the bracket is \(\partial_\theta(\theta\det D_\theta)\). All these
expressions extend to the rank-one endpoint \(D_0\), because the
determinant is a polynomial. In the same spirit, \eqref{eq:energy-J}
is the counterpart of Aubin's \(J\)-functional~\cite[p.~145]{aubin}, which
subtracts the segment energy from its linearization at the reference
potential. In our case the reference is \(B_0+u\) and the other endpoint
is \(B_0\).
\end{remark}

\begin{lemma}\label{lem:energy-balance}
The energy \(F\) is smooth and satisfies
\begin{equation}\label{eq:energy-balance}
 \Delta_z F+\int_{\Omega_z} A:\Sigma\,dx\,dP\,dZ
 =\int \frac{u}{r}\,G(r)^* A_w G(r)\,\delta(r)\,dx\,dP\,dZ .
\end{equation}
Here \(u/r\) denotes its smooth extension to the wall. The gradient
\(G(r)\) is the full gradient. In the integral on the right,
\(\delta(r)\,dx\,dP\,dZ\) is the surface measure on the fixed
fiber wall divided by \(|\nabla_{\mathrm{fib}}r|\). Moreover,
\begin{equation}\label{eq:energy-comparison}
 \frac{n-1}{n+1}F_0\leq F\leq\frac{n}{n+1}F_0 .
\end{equation}
\end{lemma}

\begin{proof}
\emph{Integration by parts on the moving fibers.}
Let \(\Omega\) be the total domain, with coordinates
\((z,x,P,Z)\) satisfying \((x,P,Z)\in\Omega_z\).
Test in the base against a smooth compactly supported \(\eta(z)\).
All boundary terms come from the wall \(r=0\), since the test function
vanishes near the ends of the base domain. If \(\mathcal A\)
is a Hermitian matrix field with divergence-free real columns, apply
\eqref{eq:green} on \(\Omega\) with \(v=\eta\).
On the wall the term with \(uG(\eta)\) vanishes because \(u=0\).
Writing \(\mathbf n=\nabla r/|\nabla r|\) for the outward real
normal and \(\nu=G(r)/|\nabla r|\) for its quaternionic column, we
also have \(G(u)=\nu\,\partial_{\mathbf n}u\) on the wall.
Rearranging Green's identity therefore gives
\begin{equation}\label{eq:energy-green}
 \int_{\Omega} u\mathcal A_{tt}\Delta_z\eta
 =\int_{\Omega} \eta\,\mathcal A:H(u)
 -\int_{r=0}\eta\,(\nu^*\mathcal A\nu)\,\partial_{\mathbf n}u\,dS .
\end{equation}
Volume measure on \(\Omega\) is \(dz\,dx\,dP\,dZ\) and \(dS\) is
the surface measure of the wall in the total domain.

To turn the boundary term into an integral over the fiber walls, first
write \(u=r(u/r)\) near \(r=0\). The quotient is smooth by taking
\(r\) as a local normal coordinate and integrating the derivative of \(u\)
along that coordinate from \(0\) to \(r\). Differentiation gives
\(\partial_{\mathbf n}u=(u/r)|\nabla r|\) at the wall. In the total
domain the coarea formula reads
\[
 dS=|\nabla r|\,\delta(r)\,dz\,dx\,dP\,dZ .
\]
Consequently the two factors of \(|\nabla r|\), one from the normal
derivative and one from the measure, cancel the denominator in
\(\nu^*\mathcal A\nu\). The remaining density in the boundary integral is
\[
 \eta\,\frac{u}{r}\,G(r)^*\mathcal A G(r)\,\delta(r)\,dz\,dx\,dP\,dZ .
\]
With \(z\) fixed, apply coarea in just the fiber variables to get
\[
 \delta(r)\,dx\,dP\,dZ
 =\frac{dS_{\mathrm{fib}}}{|\nabla_{\mathrm{fib}}r|}.
\]
This step accounts for the two different normals in
\eqref{eq:energy-balance}. The current contracts against a normal to
the wall in the total domain, while the integral is taken on one fiber.

\step{The two Laplacians.}
Take \(\mathcal A=A\) in \eqref{eq:energy-green}. The left-hand
side is \(-\int F_0(z)\Delta_z\eta(z)\,dz\). After multiplying by
\(-1\), the volume term on the right has a minus sign and the wall
term a plus sign. The cofactor regularization annihilates \(G(r)\)
by Lemma~\ref{lem:fiber-geometry}, so only \(A_w\) appears at the
wall. Using \eqref{eq:energy-schur}, we obtain the distributional identity
\[
 \begin{aligned}
 \Delta_z F_0={}&
 \int \frac{u}{r}\,G(r)^*A_w G(r)\,\delta(r)\,dx\,dP\,dZ\\
 &-\int_{\Omega_z} A:\Sigma\,dx\,dP\,dZ
   -\lambda^{-1}\int_{\Omega_z}\det H(B_0+u)\,dx\,dP\,dZ .
 \end{aligned}
\]

For \(E\), use \(\mathcal A=C(H(B_0+\theta u))\)
in \eqref{eq:energy-green}. Its real columns are divergence free by
Lemma~\ref{lem:cofactor-divergence}. (For fixed product domains in the
complex setting, Åhag et al.\ relate the Laplacian of a slice energy to
the full Monge--Ampère measure \cite[Theorem~3.1]{ahag-et-al}.) Here the
domain moves with \(z\), so we must check that the cofactor contributes
no wall term.

The coefficient \(\nu^*\mathcal A\nu\) is the determinant of the form
\(H(B_0+\theta u)\) compressed to the quaternionic orthogonal
complement \(\nu^\perp\), by Lemma~\ref{lem:cofactor}. For a vector in
\(\nu^\perp\), the product rule terms coming from a derivative of
\(r\) paired with a derivative of \(u/r\) vanish in the compression.
In detail, \eqref{eq:boundary-compression} gives
\[
 H(B_0+\theta u)|_{\nu^\perp}
   = \bigl(H(B_0)+\theta(u/r)H(r)\bigr)|_{\nu^\perp}.
\]
By Lemma~\ref{lem:fiber-geometry} the null line of \(H(r)\) lies in
\(\nu^\perp\) and is also in the kernel of \(H(|\Phi|^2)\).
Since \(B_0=\lambda|\Phi|^2\), every vector on that line is
annihilated by both terms of the compressed form. The determinant
of the compression is therefore zero. This argument uses only the
null line and applies even when the full Hessian is indefinite.

For the interior terms, the same cofactor restriction formula with
\(e_t\) in place of \(\nu\) gives \(\mathcal A_{tt}=\det D_\theta\).
The derivative formula for the cofactor gives
\(\mathcal A:H(u)=\partial_\theta\det H(B_0+\theta u)\).
Integrating \eqref{eq:energy-green} in \(\theta\) and dividing by
\(\lambda\) now yields
\[
 \begin{aligned}
 \Delta_z E
 &=\lambda^{-1}\int_{\Omega_z}
       \int_0^1\partial_\theta\det H(B_0+\theta u)\,d\theta\,dx\,dP\,dZ\\
 &=\lambda^{-1}\int_{\Omega_z}\det H(B_0+u)\,dx\,dP\,dZ.
 \end{aligned}
\]
We used \(\det H(B_0)=0\) in the second equality. The full determinant
now enters \(\Delta_z E\) and \(\Delta_z F_0\) with opposite
coefficients, \(\lambda^{-1}\) and \(-\lambda^{-1}\), so it cancels
in their sum. This cancellation follows from the first variation of the
determinant and holds even for a single Fourier mode.

The inverse chart in
\((W,Q)\) reduces the fiber integrals to a fixed domain with smoothly
varying integrands and Jacobians. The same holds for wall integrals
by using the fixed sphere in that chart. Together with smooth
dependence of \(u\) up to the wall, this shows that the distributional
identity holds pointwise. In particular,
\(A:\Sigma\geq0\) gives the wall flux as an upper bound for
\(\Delta_z F\).

\step{Size of the correction.}
The last issue is whether subtracting the determinant energy can
lose too much of \(F_0\). On a fixed fiber,
\[
 D_\theta=\theta D_1+(1-\theta)H_{\mathrm{fib}}(B_0).
\]
Since the second term has rank one, adding it changes the determinant
of \(\theta D_1\) by just its first variation
in that direction. The cofactor of \(\theta D_1\) scales as
\(\theta^{n-1}\), so
\[
 \det D_\theta
 =\theta^n\det D_1+(1-\theta)\theta^{n-1}C(D_1):H_{\mathrm{fib}}(B_0).
\]
The formula extends to \(\theta=0\) by polynomiality. Integrating
in \(\theta\), using \(\det D_1=\lambda\mu\) and
\(F=F_0+E\), gives
\[
 F=\frac{n}{n+1}F_0
 -\frac{1}{\lambda n(n+1)}
    \int_{\Omega_z}(-u)C(D_1):H_{\mathrm{fib}}(B_0)\,dx\,dP\,dZ.
\]
Here the integrals of \(\theta^n\) and
\((1-\theta)\theta^{n-1}\) are \(1/(n+1)\) and \(1/[n(n+1)]\),
respectively. Both matrices in the last contraction are positive
semidefinite, and \(-u\geq0\), giving the upper bound.

For the lower bound, use the Euler identity
\(C(D_1):D_1=n\det D_1\), obtained by differentiating
\(\det(\theta D_1)=\theta^n\det D_1\) at \(\theta=1\).
Subtracting the background contribution from the Euler identity
and integrating gives
\[
 \begin{aligned}
 &n\lambda F_0-\int_{\Omega_z}(-u)C(D_1):H_{\mathrm{fib}}(B_0)\,dx\,dP\,dZ\\
 &\qquad=\int_{\Omega_z}(-u)C(D_1):H_{\mathrm{fib}}(u)\,dx\,dP\,dZ\\
 &\qquad=\int_{\Omega_z}G_{\mathrm{fib}}(u)^*C(D_1)G_{\mathrm{fib}}(u)
          \,dx\,dP\,dZ\geq0 .
 \end{aligned}
\]
In the second equality we integrated by parts on the fiber, using
the vanishing real column divergence of \(C(D_1)\) there and the
zero boundary value of \(u\). The resulting integrand is
nonnegative because \(D_1>0\) implies \(C(D_1)>0\).
The last two displays give the exact comparison
\[
 F=\frac{n-1}{n+1}F_0
    +\frac{1}{\lambda n(n+1)}
      \int_{\Omega_z}G_{\mathrm{fib}}(u)^*C(D_1)G_{\mathrm{fib}}(u)
           \,dx\,dP\,dZ .
\]
This proves \eqref{eq:energy-comparison}.
\end{proof}

\subsection{Estimating one fiber}\label{sec:one-fiber}

We need two kinds of estimates from the same filter. Its small value
on the wall makes the boundary term in \eqref{eq:energy-balance} small.
Its reciprocal near the center lets us recover the original wave from
the current. Both estimates use a spatial kernel at scale \(1/R\).
In the upper bound we then average at the larger scale \(Ms\), where \(M\)
is a large constant, so that the derivatives of the comparison potential
are small.

\begin{proposition}\label{prop:fiber-estimates}
For a sufficiently small fixed cap aperture \(\gamma\) and sufficiently large fixed \(K\), the energy constructed above satisfies
\[
 \begin{gathered}
 F(0,\tau)\lesssim a^2\|f\|_{p_c}^{p_c}
 \qquad (|\tau|=4a),\\
 F(t,\tau)\gtrsim a^2\|e^{-it\sqrt{-\Delta}}f\|_{p_c}^{p_c}
 \qquad (|\tau|=2a),\\
 \Delta_z F(t,\tau)\leq Cs^{-2}e^{-cRs}F(t,\tau)
 \qquad (2a\leq s=|\tau|\leq3).
 \end{gathered}
\]
The estimates are uniform in \(R\), the direction of the cap, and \(t\in\mathbb R\).
\end{proposition}

\begin{proof}
Throughout the proof \(2a\leq s\leq3\), and \(N\) denotes a fixed nonnegative exponent that may increase from line to line.

\step{Upper bounds for the depth and the flux.}
For a large constant \(M\), take the positive spatial average
\[
 \Pi(x)=\big(\mathcal K_{Ms}*_x |e^{-it\sqrt{-\Delta}}f|^2(x)\big)^{1/d}.
\]
Derivatives of \(\mathcal K_{Ms}\) are bounded pointwise by the kernel times the corresponding power of \((Ms)^{-1}\). The chain rule and Young's inequality give
\begin{equation}\label{eq:depth-average}
 |\nabla_x\Pi|\lesssim (Ms)^{-1}\Pi,\qquad
 |\nabla_x^2\Pi|\lesssim (Ms)^{-2}\Pi,\qquad
 \int \Pi^{d+1}\,dx\leq\|e^{-it\sqrt{-\Delta}}f\|_{p_c}^{p_c}.
\end{equation}
For the last inequality, the exponent on the convolution is \((d+1)/d=p_c/2\).

In the wave current every scalar multiplier has a bounded polarization factor times
\[
 (\kappa/|Y|)^{d/2}\exp\!\left(-\frac{\kappa|Y|}{2}|W-v|^2\right).
\]
Cut off the symbol smoothly in \(\xi/R\), with the cutoff equal to one on the data and supported in a fixed enlarged cap \(|v|<1/30\). Up to order \(2n+4\), the symbol derivatives in \(\xi/R\) are bounded by \(C(R/s)^{d/2}(Rs)^N\), since \(\kappa\asymp R\) and \(|Y|\asymp s\). On the part of the wall where \(\chi_0(Q)\) can be nonzero, \(|Q|\leq1/20\) forces \(|W|\geq\sqrt{3}/20\). The enlarged cap still satisfies
\( |W-v|\geq\sqrt{3}/20-1/30>0\) there. Differentiating the symbol
brings out powers of \(Rs\) while leaving the Gaussian in place. This
gives the additional factor \(e^{-c_{\mathrm{wall}}Rs}\), with
\(c_{\mathrm{wall}}>0\) independent of \(\gamma\).

We use the following kernel estimate for both this multiplier and its reciprocal below. A compactly supported symbol in \(\xi/R\) whose derivatives through order \(2n+4\) are bounded by \(B\) has kernel dominated by \(CB\mathcal K_{1/R}\). This follows by applying \((1-\Delta)^{n+2}\) with respect to \(\xi/R\) in Fourier inversion and then periodizing. Cauchy--Schwarz against the probability kernel bounds the square of its output by \(CB^2 \mathcal K_{1/R}*_x|e^{-it\sqrt{-\Delta}}f|^2\). The extra convolution in \(A_w\) is harmless: \(\mathcal K_{1/R}*\mathcal K_{1/R}\lesssim \mathcal K_{1/R}\) follows by splitting the Euclidean convolution where one of the two distances to the origin is at least half the output distance, and the bound persists on periodizing. Also \(\mathcal K_{1/R}\leq (MRs)^n \mathcal K_{Ms}\). Applying these facts componentwise gives
\begin{equation}\label{eq:wave-wall-bound}
 \begin{aligned}
 \|A_w\|&\lesssim s^{-2d}(Rs)^N \Pi^d && (r\leq0),\\
 \|A_w\|&\lesssim s^{-2d}(Rs)^N e^{-c_{\mathrm{wall}}Rs}\Pi^d &&(r=0).
 \end{aligned}
\end{equation}
We have absorbed fixed powers of \(Rs\) and decreased \(c_{\mathrm{wall}}\) if necessary. The first bound also applies to \(\mu\). Its additional term is \(O(R^{-50n}a^{-2}s^{-2d})\), while the nearest lattice term in \(\mathcal K_{Ms}\) is bounded below by \(c_Ms^{n+4}\), and so \(\Pi^d\gtrsim_M R^{-n-4}\) because \(s\geq 1/R\) and \(\|e^{-it\sqrt{-\Delta}}f\|_2=1\).

The comparison function uses the same average:
\[
 H_{\mathrm{fib}}(r\Pi)\geq(1-C/M)\Pi H_{\mathrm{fib}}(r).
\]
To check this, evaluate the Hessian product rule on a quaternionic column \(w\). By Lemma~\ref{lem:fiber-geometry}, the pairing of \(G_{\mathrm{fib}}(r)\) with \(w\) is bounded by \((w^*H_{\mathrm{fib}}(r)w)^{1/2}\). By \eqref{eq:depth-average}, the pairing of \(G_{\mathrm{fib}}(\Pi)\) with \(w\) is bounded by \(CM^{-1}\Pi\) times the same expression. This bounds the two mixed terms, and the term \(rH_{\mathrm{fib}}(\Pi)\) has relative size \(O(M^{-2})\), since \(|r|\leq1\).

Fix \(M\) large enough that \(1-C/M\geq1/2\). After adding the rank-one background, the determinant of the Hessian of \(B_0+C_1(Rs)^{N_1}r\Pi\) is at least
\[
 c\lambda\big(C_1(Rs)^{N_1} \Pi s^{-2}\big)^d\geq\lambda\mu
\]
for sufficiently large fixed \(C_1\) and \(N_1\). The first inequality uses the determinant computation \eqref{eq:background-determinant} together with an isotropic lower bound for the comparison Hessian. The two potentials have the same boundary value \(B_0\). The larger determinant puts the comparison function below \(u\) by Lemma~\ref{lem:comparison}, so division by \(r<0\) and continuity to the wall yield
\[
 0\leq u/r\lesssim (Rs)^N\Pi .
\]
Consequently \(-u\lesssim (Rs)^N \Pi\). Insert this and the bound for \(\mu\) from \eqref{eq:wave-wall-bound} into \eqref{eq:energy-comparison}. The \((P,Z)\)-volume of \(\Omega_z\) is \(O(a^2s^{2d})\), which cancels the factor \(s^{-2d}\) in the density estimate. Using \eqref{eq:depth-average} gives
\begin{equation}\label{eq:fiber-upper}
 F(t,\tau)\lesssim a^2(Rs)^N\|e^{-it\sqrt{-\Delta}}f\|_{p_c}^{p_c}.
\end{equation}

For the wall integral in \eqref{eq:energy-balance} only \(G(100|W|^2)\) contributes to \(G(r)^*A_wG(r)\), because the current annihilates the gradient of the passive part \(100|Q|^2\). The full gradient of \(100|W|^2\) is \(O(s^{-1})\). Coarea in the fixed \((W,Q)\)-ellipsoid shows that the \(\delta(r)\,dP\,dZ\) mass is \(O(a^2s^{2d})\). The wall bound in \eqref{eq:wave-wall-bound}, the depth comparison, and \eqref{eq:depth-average} now imply
\begin{equation}\label{eq:fiber-flux}
 \Delta_z F(t,\tau)
 \lesssim a^2s^{-2}(Rs)^N e^{-c_{\mathrm{wall}}Rs}\|e^{-it\sqrt{-\Delta}}f\|_{p_c}^{p_c}.
\end{equation}

\step{Recovering the wave in the interior.}
Take a center \((x_0,P_0,Z_0)\) with \(|W|\leq\gamma\) and \(|Q|\leq1/400\). On the ball of radius \(1/R\) in the fiber centered there we have \(|W|\leq2\gamma\) and \(|Q|\leq1/100\) once \(K\gamma\) and \(K\) are large, because \(W\) changes by \(O((Rs)^{-1})\) and \(\Phi\) by \(O(1/R)\) on such a ball. In particular \(\chi_0(Q)=1\) there. In the first multiplier \(S_0\) from \eqref{eq:time-component}, the factor
\[
 b_0=\frac{1+\operatorname{Re}(W)\cdot v/2}{\sqrt{2}}
\]
is bounded away from zero. Multiply the reciprocal symbol by a smooth cutoff equal to one on the data and supported on \(|v|\leq2\gamma\) in a fixed enlarged annulus. The derivatives of the resulting symbol through order \(2n+4\) are bounded by
\[
 C_\gamma(s/R)^{d/2}(Rs)^N e^{C\gamma^2Rs}.
\]
The exponential has a constant independent of \(\gamma\), since \(|W-v|\leq4\gamma\); the cutoff derivatives affect only the prefactor. At any fixed pair \((P,Z)\) in the ball, the inverse multiplier
composed with \(S_0\) is the identity on the frequency support.
The kernel bound and Cauchy--Schwarz applied at \(x_0\) give
\[
 |e^{-it\sqrt{-\Delta}}f(x_0)|^2
 \leq C_\gamma(s/R)^d(Rs)^N e^{C\gamma^2Rs}
   \big(\mathcal K_{1/R}*_x
       |S_0 e^{-it\sqrt{-\Delta}}f|^2\big)(x_0).
\]
For \(|x-x_0|\leq1/R\), the positive averaging kernel at \(x-y\)
is comparable to its value at \(x_0-y\), uniformly in \(y\):
these shifts correspond to a bounded translation in the
rescaled Euclidean kernel, and the comparison holds term by term
on periodizing. We can therefore replace evaluation at \(x_0\) by
evaluation at \(x\) on the right. By \eqref{eq:time-component} and
\(\chi_0(Q)=1\), the resulting average is bounded by \(\mu(x,P,Z)\).
The wave amplitude at the center thus bounds the density from below
on the whole ball: \(\mu\geq\underline\mu(x_0)\) there, where
\[
 \underline\mu(x_0)=c_\gamma(R/s)^d(Rs)^{-N}e^{-C\gamma^2Rs}
       |e^{-it\sqrt{-\Delta}}f(x_0)|^2 .
\]

If \(\underline\mu(x_0)>0\), compare with the quadratic
\[
 c_1\underline\mu(x_0)^{1/d}\bigl(|x-x_0|^2+|P-P_0|^2+|Z-Z_0|^2-R^{-2}\bigr).
\]
The distances are computed in the local real coordinates of the fiber. The Hessian of this quadratic equals \(6c_1\underline\mu(x_0)^{1/d}I\). Since \(\underline\mu(x_0)^{1/d}\leq\mu^{1/d}\leq\lambda\), the determinant after adding the background is at most \(\lambda\underline\mu(x_0)\) for a small fixed \(c_1>0\), by \eqref{eq:background-determinant}. The comparison quadratic vanishes on the boundary of the ball, where \(u\leq0\), and the solution has the larger determinant by \eqref{eq:fiber-dirichlet}. Hence \(u\) lies below the quadratic by Lemma~\ref{lem:comparison}, and at the center
\[
 (-u)\mu\geq c_1R^{-2}\underline\mu(x_0)^{1+1/d}.
\]
This also holds when \(\underline\mu(x_0)=0\). This is the point at which the exponent \(p_c\) enters: the squared wave amplitude appears to the power \(1+1/d=p_c/2\).

The \((P,Z)\)-volume of the region of centers is comparable to \(a^2s^{2d}\gamma^{2d}\). Integrating the last bound and using \eqref{eq:energy-comparison} gives
\begin{equation}\label{eq:fiber-lower}
 F(t,\tau)\gtrsim_\gamma
 a^2(Rs)^{-N}e^{-C\gamma^2Rs}\|e^{-it\sqrt{-\Delta}}f\|_{p_c}^{p_c}.
\end{equation}
Here the scale factor \(R^{-2}s^{2d}(R/s)^n=(Rs)^{n-2}\) is absorbed in the polynomial factor.

Choose \(\gamma\) small enough that \(C\gamma^2<c_{\mathrm{wall}}/2\),
then choose \(K\) large enough for the balls of radius \(1/R\) to stay in
\(|W|\leq2\gamma\) and \(|Q|\leq1/100\). With these constants fixed,
take \(R\) large enough for the density ceiling \eqref{eq:density-ceiling}
and for \(4a<1\).
Dividing \eqref{eq:fiber-flux} by \eqref{eq:fiber-lower} leaves \(s^{-2}\) times a polynomial times a decaying exponential in \(Rs\). Absorbing the polynomial into the exponential gives the differential inequality of the proposition. At \(s=4a\) and \(s=2a\), \(Rs\) takes the fixed values \(4K\) and \(2K\), so the two circle estimates follow from \eqref{eq:fiber-upper} and \eqref{eq:fiber-lower}.
\end{proof}

\subsection{A scalar comparison in the base variables}\label{sec:cylinder}

The differential inequality for the energy is elliptic in the three base
variables \((t,\tau_1,\tau_2)\). In this part of the proof, real time is
one of the coordinates of a scalar elliptic problem. Since the potential
\(V\) is nonnegative, the inequality \(\Delta_zF\leq VF\) makes \(F\) a
supersolution of \(-\Delta_z+V\), and the maximum principle bounds \(F\)
from below by any solution of the homogeneous equation that lies below
\(F\) on the boundary. We use this direction of comparison
to bound the energy at \(t=0\) from below by the energy on \(1\leq t\leq2\).

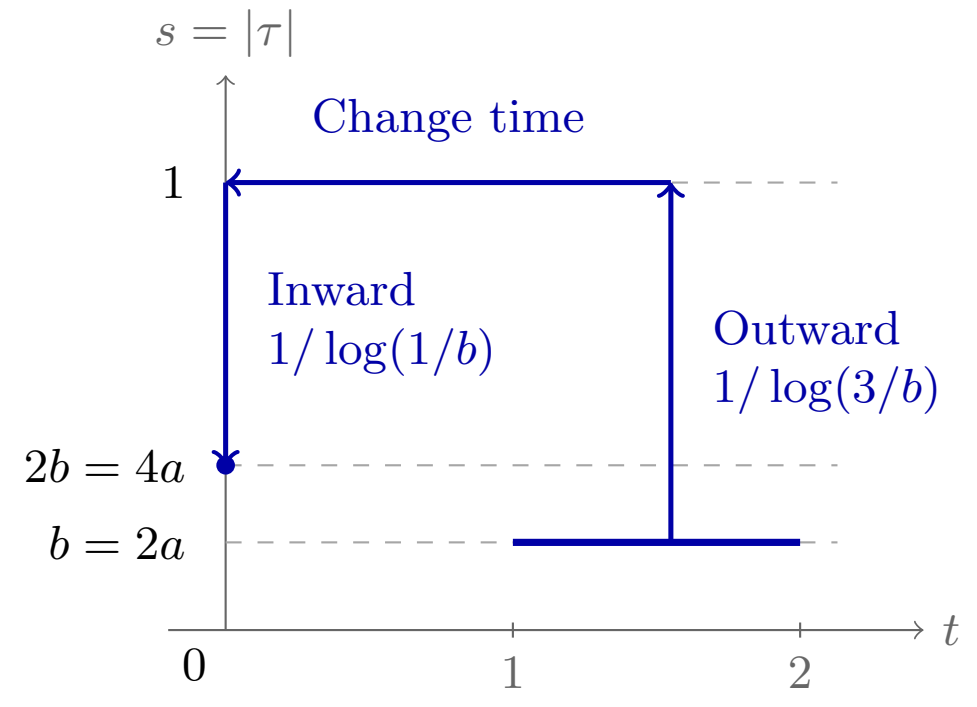
\begin{figure}[tbp]
\centering
\begin{tikzpicture}[x=1.9cm,y=1cm,font=\small]
  \draw[->,black!60] (-.20,.2) -- (2.43,.2) node[right] {\(t\)};
  \draw[->,black!60] (0,.2) -- (0,3.87) node[above] {\(s=|\tau|\)};
  \foreach \x in {1,2}
    \draw[black!60] (\x,.15) -- (\x,.25) node[below=3pt] {\(\x\)};
  \node[below left] at (0,.2) {\(0\)};
  \draw[dashed,black!35] (0,.78) -- (2.13,.78);
  \draw[dashed,black!35] (0,1.29) -- (2.13,1.29);
  \draw[dashed,black!35] (0,3.16) -- (2.13,3.16);
  \node[left=4pt] at (0,.78) {\(b=2a\)};
  \node[left=4pt] at (0,1.29) {\(2b=4a\)};
  \node[left=4pt] at (0,3.16) {\(1\)};
  \draw[blue!65!black,line width=1.35pt] (1,.78) -- (2,.78);
  \draw[->,blue!65!black,line width=.95pt] (1.55,.78) -- (1.55,3.16)
     node[midway,right=4pt,align=left] {Outward\\ \(1/\log(3/b)\)};
  \draw[->,blue!65!black,line width=.95pt] (1.55,3.16) -- (0,3.16)
     node[midway,above=4pt] {Change time};
  \draw[->,blue!65!black,line width=.95pt] (0,3.16) -- (0,1.29)
     node[midway,right=4pt,align=left] {Inward\\ \(1/\log(1/b)\)};
  \fill[blue!65!black] (0,1.29) circle[radius=1.75pt];
\end{tikzpicture}
\caption{The comparison in \((t,s)\) after averaging over the angle of
\(\tau\). The starting segment represents \(1\leq t\leq2\) on the circle
\(s=b\). The endpoint represents the circle \(s=2b\) at \(t=0\). The two
radial comparisons supply one logarithmic factor each. The transfer in
time uses a local mean value inequality at fixed scale. The radii are
schematic.}
\label{fig:cylinder}
\end{figure}

The energy estimates in Proposition~\ref{prop:fiber-estimates} involve two different circles in the \(\tau\)-plane. On the smaller circle the energy controls the wave for \(1\le t\le 2\). On the larger circle the energy at \(t=0\) is bounded by the initial data. To connect the two, we first move out to radius \(1\), where the geometry is independent of \(R\), and use a local mean value inequality to move backward in real time. We then move inward to reach the larger of the two small circles, as
shown in Figure~\ref{fig:cylinder}.

Readers who know the induction-on-scales arguments of Wolff~\cite{wolff} and of Bourgain and Demeter~\cite{bdec} may find the following heuristic comparison helpful. At auxiliary radius \(s\) the filter resolves angles at scale \((Rs)^{-1/2}\). On the circles \(s=b=2a\) and \(s=2b\) this scale does not depend on \(R\), and Proposition~\ref{prop:fiber-estimates} compares the energy there with \(L^{p_c}\) norms of the wave. At \(s=1\) it is the wave packet scale \(R^{-1/2}\), and there the proof moves in time without a logarithmic loss. An induction on scales passes from one scale to the next and pays a price at each step, which is why such arguments end with losses of the form \(R^\varepsilon\). Here each radial leg passes through all intermediate angular scales at once, by a single comparison for the scalar inequality \(\Delta_zF\leq VF\), and costs one factor of \(\log(1/b)\asymp\log R\).

The comparison with a potential is OpenAI's cylinder estimate
\cite[Lemma~6.1]{openai}.

\begin{lemma}[Comparison on a cylinder]
\label{lem:cylinder}
Let \(0<b<1/4\) and let \(V\colon[b,3]\to[0,\infty)\) be continuous, with
\[
 \int_b^3 s\log(s/b)V(s)\,ds\le M_V,
 \qquad
 \sup_{1/2\le s\le 3/2} V(s)\le M_V.
\]
Suppose \(F\ge 0\) is smooth on \(\mathbb R_t\times\{\tau\in\mathbb R^2:b\le |\tau|\le 3\}\) and obeys
\(\Delta_z F\le V(|\tau|)F\) in the interior. Then its circular mean
\[
 \overline F(t,s)=\frac{1}{2\pi}\int_0^{2\pi} F(t,s\cos\theta,s\sin\theta)\,d\theta
\]
satisfies
\begin{equation}
\label{eq:cylinder-propagation}
 \overline F(0,s_0)\ge c_{M_V}\frac{\log(s_0/b)}{[\log(3/b)]^2}
      \int_1^2 \overline F(t,b)\,dt
 \qquad (b\le s_0\le 1).
\end{equation}
\end{lemma}

\begin{proof}
In polar coordinates in the \(\tau\)-plane, the angular part
of the Laplacian integrates to zero. Averaging therefore gives
\[
  \partial_t^2\overline F+\partial_s^2\overline F
        +s^{-1}\partial_s\overline F\leq V(s)\overline F .
\]

\step{The radial equation.}
For \(S=1\) or \(3\) and fixed \(\omega>0\), consider the boundary
value problem
\[
  U''+s^{-1}U'=(V(s)+\omega^2)U\qquad(b<s<S)
\]
with endpoint values \(0,1\) in either order. The term \(\omega^2\)
comes from the sine or cosine in the real time variable below.
For the Laplacian alone, the equation is \((sU')'=0\), so the
solutions with boundary values \((0,1)\) and \((1,0)\) are
\[
  \frac{\log(s/b)}{\log(S/b)}
  \qquad\text{and}\qquad
  \frac{\log(S/s)}{\log(S/b)},
\]
respectively. These are the usual harmonic functions on a planar
annulus~\cite[Theorem~3.18]{morters-peres}. We show that the nonnegative
potential decreases each function by at most a constant factor,
uniformly in the length \(\log(S/b)\) of the interval in logarithmic
coordinates.

Set \(x=\log(s/b)\) and \(D=\log(S/b)\). In these coordinates
\(U''(s)+s^{-1}U'(s)=s^{-2}\partial_x^2U\), so the coefficient of
\(U\) in the transformed equation is
\(q(x)=s^2(V(s)+\omega^2)\), where \(s=be^x\).
Solve the initial value problem \(\psi''=q\psi\),
\(\psi(0)=0\), \(\psi'(0)=1\). Since \(q\geq0\), \(\psi\) is
convex as long as it is nonnegative, so \(\psi(x)\geq x\) throughout
\([0,D]\). Integrating the equation twice gives
\[
 \frac{\psi(x)}{x}
  =1+\int_0^x \frac{y(x-y)}{x}q(y)\frac{\psi(y)}{y}\,dy.
\]
For \(0<y<x\leq D\) the factor \(y(x-y)/x\) is at most
\(\min(y,D-y)\). The hypothesis controls
\(\int_0^D \min(y,D-y) s^2V(s)\,dy\) by \(M_V\), since
\(\min(y,D-y)\leq y=\log(s/b)\) and \(dy=ds/s\).
For the constant \(\omega^2\) term, use
\(s=Se^{-(D-y)}\) and integrate
\(\omega^2 S^2(D-y)e^{-2(D-y)}\) to bound its contribution by
\(\omega^2S^2/4\).
Gronwall's inequality gives \(\psi(D)/D\leq C_{M_V,\omega}\), so
the solution \(\psi(x)/\psi(D)\) is bounded below by
\(c_{M_V,\omega}x/D\). For the reversed boundary values, reverse the
coordinate \(x\mapsto D-x\). The integral bound for the coefficient
is unchanged because \(\min(y,D-y)\) is symmetric. This gives
the lower bound \(c_{M_V,\omega}(D-x)/D\).

\step{Moving outward.}
Use the first Dirichlet eigenfunction of the interval \([0,3]\)
to take a weighted average in time,
\[
 \mathcal F(s)=\int_0^3\sin(\pi t/3)\overline F(t,s)\,dt .
\]
This weight is bounded away from zero on \([1,2]\).
Two integrations by parts in \(t\) give
\[
 \int_0^3\sin(\pi t/3)\partial_t^2\overline F(t,s)\,dt
 =-\frac{\pi^2}{9}\mathcal F(s)
  +\frac{\pi}{3}\bigl(\overline F(0,s)+\overline F(3,s)\bigr).
\]
The boundary contributions are nonnegative. The averaged differential
inequality thus gives \(\mathcal F''+s^{-1}\mathcal F'\leq (V+\pi^2/9)\mathcal F\).
At \(s=b\) match \(\mathcal F(b)\), and at \(s=3\) use \(\mathcal F(3)\geq0\).
Comparison with \(\mathcal F(b)\) times the radial solution with boundary
values \((1,0)\), for \(S=3\) and \(\omega=\pi/3\), gives
\[
 \mathcal F(s)\geq c_{M_V}\mathcal F(b)\frac{\log(3/s)}{\log(3/b)}.
\]
For \(3/4\leq s\leq5/4\), the numerator is comparable to one.
Integrate over this annulus with the polar measure \(2\pi s\,ds\).
Since the sine is at most one, the integral of \(\overline F(t,|\tau|)\)
over \([0,3]\times\{3/4\leq|\tau|\leq5/4\}\) is bounded below by
\(c_{M_V}\mathcal F(b)/\log(3/b)\). This gives the first logarithmic factor.

\step{Changing time at radius one.}
In \(-2<t<4\), \(1/2<|\tau|<3/2\), we have
\(\Delta_z\overline F\leq M_V\overline F\).
For a local inequality with a constant independent of \(b\), add
a real coordinate \(\sigma\). In a small slab about \(\sigma=0\), the function
\(\overline F(t,|\tau|)\cos(\sqrt{M_V+1}\,\sigma)\) is nonnegative and its
four-dimensional Laplacian is nonpositive, because the second
derivative in \(\sigma\) contributes \(-(M_V+1)\) times the function.
By the mean value inequality for superharmonic functions, its value at
\((t,\tau,0)\) bounds from below its average over a ball centered there.
On a product of a smaller
three-dimensional ball and a short interval in \(\sigma\) inside that ball,
the cosine is bounded away from zero. Thus \(\overline F(t,|\tau|)\)
bounds a fixed multiple of its local three-dimensional average.

Cover \([-1,3]\times\{3/4\leq|\tau|\leq5/4\}\) with finitely many
overlapping balls of one fixed radius small enough for this estimate.
At least one ball carries a fixed fraction of the integral just
obtained. For a point in a neighboring ball, a ball of four times the
covering radius about that point contains the first ball. The mean
value inequality gives a lower bound there in terms of the integral on the
first ball. Iterating along the finite cover gives
\[
 \overline F(t,1)\geq
 c_{M_V}\frac{\mathcal F(b)}{\log(3/b)}
 \qquad (|t|\leq1).
\]
The radii and the number of balls depend only on \(M_V\) and the fixed
annulus, so this step incurs no additional logarithmic factor.

\step{Moving inward.}
Take the radial solution \(U\) on \([b,1]\) with \(U(b)=0\), \(U(1)=1\),
and \(\omega=\pi/2\). The time factor \(\cos(\pi t/2)\) makes
\[
 c_{M_V}\frac{\mathcal F(b)}{\log(3/b)}
     \cos(\pi t/2)U(|\tau|)
\]
a function annihilated by \(-\Delta_z+V\) in the cylinder
\(|t|<1\), \(b<|\tau|<1\). It vanishes on \(|t|=1\) and \(|\tau|=b\),
and on \(|\tau|=1\) it lies below \(\overline F\).
The maximum principle and the radial lower bound at \((0,s_0)\) give
\[
 \overline F(0,s_0)\geq
 c_{M_V}\frac{\mathcal F(b)}{\log(3/b)}
          \frac{\log(s_0/b)}{\log(1/b)} .
\]
At \(s_0=2b\) the numerator \(\log(s_0/b)\) is a fixed positive number, so
in the application this step costs the second logarithmic factor. Since
\(\mathcal F(b)\gtrsim\int_1^2\overline F(t,b)\,dt\) and
\(\log(1/b)\leq\log(3/b)\), this proves
\eqref{eq:cylinder-propagation}.
\end{proof}

\subsection{From the determinant energy to the annular estimate}\label{sec:annular-proof}

We first work on the fixed torus \(\mathbb T^n=\mathbb R^n/(100\mathbb Z)^n\) and in one cap from Proposition~\ref{prop:fiber-estimates}. By homogeneity assume \(\|f\|_2=1\) and take \(R\) large enough that \(b=2a=2K/R<1/4\). The potential \(V(s)=Cs^{-2}e^{-cRs}\) meets the hypotheses of Lemma~\ref{lem:cylinder} with a constant \(M_V\) independent of \(R\): the substitution \(s=b\rho\) gives
\[
 \int_b^3 s\log(s/b)V(s)\,ds
 \le C\int_1^\infty \frac{\log \rho}{\rho}e^{-2cK\rho}\,d\rho ,
\]
and \(V\) is uniformly bounded on \([1/2,3/2]\). The two circle estimates remain true after averaging \(F(t,\tau)\) in the angle of \(\tau\), so the resulting circular mean satisfies
\[
 \overline F(0,2b)\lesssim a^2\|f\|_{p_c}^{p_c},
 \qquad
 \overline F(t,b)\gtrsim a^2\|e^{-it\sqrt{-\Delta}}f\|_{p_c}^{p_c}.
\]
Apply \eqref{eq:cylinder-propagation} at \(s_0=2b\) and cancel \(a^2\). After taking \(p_c\)-th roots we obtain the one-cap estimate with the factor \([\log(2+R)]^{2/p_c}\), and homogeneity removes the \(L^2\) normalization. For bounded \(R\geq2\), the frequencies lie in a fixed finite subset of the torus lattice, so the same estimate follows directly from the Fourier series.

We now sum over caps on this same torus. Write
\(\phi=\sum_\theta\phi_\theta\) using a smooth partition into caps of
fixed aperture, so the number of terms is independent of \(R\).
For any periodic \(f\), the multiplier with symbol
\(\phi_\theta(\xi/R)\) selects modes in the corresponding cap from
the fixed frequency lattice. Its Euclidean kernel is
\(R^n\check\phi_\theta(Rx)\). To get its torus kernel, sum the
translates by the period lattice \((100\mathbb Z)^n\). Integrating their
absolute values over a period cell gives at most
\(\|\check\phi_\theta\|_{L^1(\mathbb R^n)}\), since the translated
cells tile \(\mathbb R^n\). Applying the one-cap estimate to each
filtered input and summing gives
\begin{equation}
\label{eq:torus-annular}
 \|e^{-it\sqrt{-\Delta}}\phi_R f\|_{L^{p_c}([1,2]\times\mathbb T^n)}
 \lesssim [\log(2+R)]^{2/p_c}\|f\|_{L^{p_c}(\mathbb T^n)}
 \qquad (R\ge 2).
\end{equation}

\begin{proof}[Proof of Theorem~\ref{thm:annular}]
On \(\mathbb R^n\) the kernel of \(e^{-it\sqrt{-\Delta}}\phi_R\) is
\[
 \Gamma_R(t,x)=\frac{R^n}{(2\pi)^n}\int_{\mathbb R^n}
      e^{iR(x\cdot\eta-t|\eta|)}\phi(\eta)\,d\eta .
\]
For \(1\le t\le2\) and \(|x|\ge3\), the gradient of \(x\cdot\eta-t|\eta|\) in \(\eta\) has size at least \(|x|/3\). Integrating by parts repeatedly in the gradient direction, using that the support of \(\phi\) stays away from \(0\), gives
\begin{equation}
\label{eq:wave-kernel-away}
 |\Gamma_R(t,x)|\le C_N R^{-N}(1+|x|)^{-N}
 \qquad (1\le t\le 2,\ |x|\ge 3)
\end{equation}
for every \(N>0\). It suffices to integrate by parts \(N+n\) times, which offsets the prefactor \(R^n\).

Tile \(\mathbb R^n\) by unit cubes \(\mathcal Q\) and let \(10\mathcal Q\) be the concentric cube of side ten. For \(x\in \mathcal Q\), the contribution from \(f\mathbf 1_{(10\mathcal Q)^c}\) is bounded by a single function independent of \(\mathcal Q\),
\[
 C_N R^{-N}\bigl((1+|\cdot|)^{-N}*|f|\bigr)(x),
\]
because \(|x-y|\ge9/2\) for this contribution, so that \eqref{eq:wave-kernel-away} applies. Young's inequality with \(N>n\) controls the \(L^{p_c}\) spacetime norm of all the far contributions by \(O(R^{-N})\|f\|_{p_c}\).

For the nearby data, let \(x_{\mathcal Q}\) be the center of \(\mathcal Q\) and
use the translated coordinates in which \(\mathcal Q-x_{\mathcal Q}=[-1/2,1/2)^n\).
On the period cell \([-50,50)^n\) set
\[
 f_{\mathcal Q}(y)=f(x_{\mathcal Q}+y)\mathbf 1_{[-5,5)^n}(y),
\]
and extend it periodically. Its copies have disjoint supports, so with
Lebesgue measure on the period cell,
\(\|f_{\mathcal Q}\|_{L^{p_c}(\mathbb T^n)}
  =\|f\|_{L^{p_c}(10\mathcal Q)}\).
The spatial cutoff can spread the Fourier support of \(f_{\mathcal Q}\), but
estimate \eqref{eq:torus-annular} applies directly to this periodic
input because \(\phi_R\) acts after the spatial cutoff.

To compare the two outputs, periodize \(\Gamma_R(t,\cdot)\). Its periodization
has Fourier coefficient \(100^{-n}e^{-it|\xi|}\phi(\xi/R)\) at each
\(\xi\in(2\pi/100)\mathbb Z^n\), by integrating over the translated
period cells. Thus it is the torus convolution kernel with respect to
Lebesgue measure. For \(x\in[-1/2,1/2)^n\) the torus output is
\[
 (e^{-it\sqrt{-\Delta}}\phi_R f_{\mathcal Q})(x)
  =\sum_{\ell\in\mathbb Z^n}\int_{[-5,5)^n}
       \Gamma_R(t,x-y+100\ell) f(x_{\mathcal Q}+y)\,dy .
\]
The term \(\ell=0\) is exactly the Euclidean output for
\(f\mathbf 1_{10\mathcal Q}\) evaluated at \(x_{\mathcal Q}+x\). For the other
terms, every coordinate of \(x-y\) is at most \(11/2\) in absolute value.
In each coordinate in which \(\ell\) is nonzero, the absolute value of the
coordinate of \(x-y+100\ell\) is then at least \(100-11/2=94.5\)
times that of the coordinate of \(\ell\). Squaring and summing over those
coordinates gives \(|x-y+100\ell|\ge 94.5|\ell|\). By
\eqref{eq:wave-kernel-away}, the total contribution of the nonzero
translates is therefore bounded pointwise by
\[
 C_N R^{-N}\|f\|_{L^1(10\mathcal Q)}
       \sum_{\ell\in\mathbb Z^n\setminus\{0\}}(1+|\ell|)^{-N}
 \lesssim_{N,n} R^{-N}\|f\|_{L^{p_c}(10\mathcal Q)} .
\]
Here \(N>n\) makes the sum finite, and Hölder's inequality applies on
a cube of fixed volume. Since \([1,2]\times\mathcal Q\) has measure
one, the same bound holds in \(L^{p_c}\) there. Applying
\eqref{eq:torus-annular} to \(f_{\mathcal Q}\) now gives
\[
 \|e^{-it\sqrt{-\Delta}}\phi_R(f\mathbf 1_{10\mathcal Q})\|_{L^{p_c}([1,2]\times \mathcal Q)}
 \lesssim [\log(2+R)]^{2/p_c}\|f\|_{L^{p_c}(10\mathcal Q)}.
\]
Raise this bound to the power \(p_c\) and sum over the disjoint output cubes.
The enlarged input cubes \(10\mathcal Q\) have bounded overlap, so
\(\sum_{\mathcal Q}\|f\|_{L^{p_c}(10\mathcal Q)}^{p_c}
\lesssim_n\|f\|_{L^{p_c}(\mathbb R^n)}^{p_c}\).
Together with the bound for the far contributions, this completes the proof.
\end{proof}

We end this subsection with the focusing example from the first paragraph of the introduction, which tests the logarithms in Theorem~\ref{thm:annular}; compare~\cite[Proposition~8.1]{openai}. Take \(\phi\) real, radial and nonnegative, not identically zero, and let \(\widetilde \Gamma_R\) be defined as \(\Gamma_R\) with \(\phi^2\) in place of \(\phi\). For \(t_0=3/2\), the function \(f=\Gamma_R(-t_0,\cdot)\) satisfies
\[
 e^{-it\sqrt{-\Delta}}\phi_Rf=\widetilde \Gamma_R(t-t_0,\cdot),
\]
a wave that focuses at the origin at time \(t_0\). By stationary phase there is a constant \(C\) such that, for \(C/R\le|t|\le2\), the kernels \(\Gamma_R(t,\cdot)\) and \(\widetilde \Gamma_R(t,\cdot)\) have size comparable to \(R^n(R|t|)^{-d/2}\) on a neighborhood of width comparable to \(1/R\) of the sphere \(|x|=|t|\), and decay rapidly away from it. Since \(p_cd/2=n\), this gives
\[
 \|\widetilde \Gamma_R(t,\cdot)\|_{p_c}^{p_c}\asymp R^{np_c-n-1}|t|^{-1}
 \quad(C/R\le|t|\le1/2),
 \qquad
 \|f\|_{p_c}^{p_c}\asymp R^{np_c-n-1}.
\]
Now integrate over \(1\le t\le2\). Each dyadic range of \(|t-t_0|\) between \(C/R\) and \(1/2\) contributes a comparable amount, and there are about \(\log R\) such ranges, so for large \(R\)
\[
 \|e^{-it\sqrt{-\Delta}}\phi_Rf\|_{L^{p_c}([1,2]\times\mathbb R^n)}^{p_c}
 \gtrsim\log R\,\|f\|_{p_c}^{p_c}.
\]
Thus the \(p_c\)-th power of the constant in Theorem~\ref{thm:annular} grows at least like \(\log R\) for this \(\phi\), whereas the proof gives \([\log(2+R)]^2\) for this power. The example shows that at least one logarithm is necessary. It does not decide whether the second one is.

\subsection{Interpolation and summation over frequencies}\label{sec:interpolation}

For the same multiplier \(e^{-it\sqrt{-\Delta}}\phi_R\), Plancherel gives the \(L^2(\mathbb R^n)\to L^2([1,2]\times\mathbb R^n)\) bound uniformly in \(R\). The \(L^\infty\to L^\infty\) bound is \(O(R^{d/2})\). To check it, split the cutoff into \(O(R^{d/2})\) pieces supported in caps of aperture comparable to \(R^{-1/2}\). In coordinates \(\eta=(R^{-1/2}\,\eta',\rho)\) centered on the axis of each cap, the rescaled cutoff is supported in a fixed compact set with \(\rho\) bounded away from \(0\) and with uniformly bounded derivatives. After removing the linear term \(-Rt\rho\), the part of the phase from the wave multiplier is
\[
 -Rt\bigl(\sqrt{\rho^2+R^{-1}|\eta'|^2}-\rho\bigr)
 =-\frac{t|\eta'|^2}{\sqrt{\rho^2+R^{-1}|\eta'|^2}+\rho}.
\]
All its derivatives are uniformly bounded on the support for \(1\le t\le2\). Integration by parts in the rescaled Fourier integral therefore bounds each kernel's \(L^1\) norm uniformly. Undoing a linear scaling of frequency or a translation of the spatial kernel preserves this norm. Summing over the caps gives the asserted bound.

Between \(2\) and \(p_c\), interpolation with Theorem~\ref{thm:annular} costs only the logarithmic factor with exponent \(d(1-2/p)\). For \(p\ge p_c\), the weight of the \(p_c\) bound in the interpolation with \(L^\infty\) is \(p_c/p\). The power of \(R\) is then
\((d/2)(1-p_c/p)=d/2-n/p\), and the power of the logarithm is \(2/p\). Thus for all \(2\le p<\infty\) we have
\begin{equation}
\label{eq:interpolated-annular}
 \begin{aligned}
 &\|e^{-it\sqrt{-\Delta}}\phi_R f\|_{L^p([1,2]\times\mathbb R^n)}\\
 &\qquad\lesssim_p R^{\sigma_p}
     [\log(2+R)]^{\min\{d(1-2/p),2/p\}}\|f\|_{L^p(\mathbb R^n)} .
 \end{aligned}
\end{equation}

\begin{proof}[Proof of Theorem~\ref{thm:local-smoothing}]
For \(2<p<\infty\), take a smooth dyadic partition of unity in spatial frequency. On the annulus \(|\xi|\asymp 2^\ell\), \(\ell\ge1\), apply \eqref{eq:interpolated-annular} with \(R=2^\ell\) and an enlarged cutoff equal to one on the annulus. The product of the annular cutoff with
\(2^{\ell s_p}(1+|\xi|^2)^{-s_p/2}\) has a uniformly integrable convolution kernel, as one sees by scaling \(\xi\) by \(2^\ell\) to a fixed compact annulus and integrating by parts. Consequently the contribution of the \(\ell\)-th piece is bounded by
\[
 C_{n,p,\varepsilon}\,2^{-\ell\varepsilon}(1+\ell)^{\min\{d(1-2/p),2/p\}}
     \|(1-\Delta)^{s_p/2}f\|_{L^p(\mathbb R^n)} .
\]
The coefficients are summable. For the low-frequency term, the symbol \(e^{-it|\xi|}\) times a smooth cutoff near the origin and \((1+|\xi|^2)^{-s_p/2}\) satisfies the Mihlin derivative bounds uniformly for \(1\le t\le2\). This controls the low frequencies by the same Sobolev norm. At \(p=2\), Plancherel gives the result directly with zero Sobolev order.
\end{proof}

\appendix

\section{Proof of \texorpdfstring{Theorem~\ref{thm:dirichlet}}{Theorem 2.5}}\label{app:dirichlet}

\begin{proof}[Proof of Theorem~\ref{thm:dirichlet}]
Uniqueness follows from Lemma~\ref{lem:comparison}. For existence we lift
the problem to four slots and use the method of continuity.

\step{The fourth slot and the continuity path.}
Adjoin one copy of \(\mathbb R/\mathbb Z\) in the \(k\)-slot of each
quaternionic coordinate. This gives the lifted domain
\(\widetilde\Omega=\Omega\times(\mathbb R/\mathbb Z)^m\) with the full
four-slot Hessian. Extend \(r\), \(B_0\) and \(g\) constantly in the added
variables. Every term involving a derivative in a fourth slot then
vanishes, so their Hessians are the original three-slot ones.
The boundary still comes entirely from \(\partial V\). With all four
slots available, quaternionic unitary rotations act in flat charts.
Their infinitesimal action will help estimate mixed second derivatives
at the boundary.

Start at the known correction \(u=r\) and vary the density along
\[
  g_\theta=(1-\theta)\det H(B_0+r)+\theta g,\qquad 0\leq\theta\leq 1.
\]
For each \(\theta\) we look for an admissible solution of
\(\log\det H(B_0+u)=\log g_\theta\) with \(u=0\) at the boundary, and we
show that the set of \(\theta\) for which one exists is open and closed
in \([0,1]\). Openness comes from the invertibility of the linearized
equation, and closedness from uniform estimates.

\step{Openness.}
At a solution \(\psi=B_0+u\), the derivative of the equation in the direction of a function \(w\) that vanishes on the boundary is
\[
  \mathcal Lw=H(\psi)^{-1}:H(w),\qquad
  \text{with real symbol }\widehat\xi^*H(\psi)^{-1}\widehat\xi .
\]
At this fixed solution the Hessian is positive on the closure
of the domain, so \(\mathcal L\) is uniformly elliptic. Use a finite cover by
flat interior charts and charts flattening the boundary in the \(V\)
factor. The local Schauder estimates~\cite[Chapter~6]{gilbarg-trudinger}
give the global estimate
\[
  \|w\|_{C^{2,\alpha}}\lesssim \|\mathcal Lw\|_{C^\alpha}+\|w\|_{C^0},
  \qquad w|_{\partial\widetilde\Omega}=0 .
\]
The constants in this estimate can be chosen uniformly along
the linear path from \(\Delta\) to \(\mathcal L\), since the ellipticity and
coefficient bounds are uniform. Also \(\Delta r=\operatorname{tr}H(r)>0\)
and \(\mathcal Lr>0\), so each operator along the path sends \(r\) to a function
bounded below by a positive constant. For a solution \(w\) with right-hand
side \(f\), apply the maximum principle to \(\pm w+C\|f\|_{C^0}r\).
This gives \(\|w\|_{C^0}\lesssim\|f\|_{C^0}\) and removes the last term
from the estimate. The starting Dirichlet Laplacian is invertible
by the variational construction followed by Schauder regularity.
The Poincaré inequality used there holds in \(V\) for each fixed value
of the torus variables. The linear continuity method now shows that
\(\mathcal L\) is an isomorphism from the \(C^{2,\alpha}\) functions that vanish on the boundary to \(C^\alpha\).
The implicit function theorem gives nearby admissible solutions,
and their smoothness follows from elliptic regularity.

\step{Barriers and interior derivatives.}
Both endpoints of the density path are strictly above \(\det H(B_0)\).
At \(\theta=0\) this follows by adding the positive definite matrix
\(H(r)\) to \(H(B_0)\geq0\). The gap is uniform on the closed domain.
For small \(b'>0\), \(\det H(B_0+b'r)\) is uniformly close to \(\det H(B_0)\).
For large \(b\), \(H(B_0+br)\geq bH(r)\) has determinant growing at least
as \(b^m\). Thus we can choose \(0<b'<1<b\) so that
\[
  \det H(B_0+b'r)< g_\theta<\det H(B_0+br)
\]
throughout the domain and along the path. All three potentials have
boundary value \(B_0\). Comparison gives \(br\leq u\leq b'r\).
Near the boundary \(u=\rho r\) for a smooth quotient \(\rho\), since
\(r\) is a local coordinate transverse to \(\partial V\). Dividing the
barriers by \(r<0\) and taking the boundary limit gives
\(b'\leq\rho\leq b\). Therefore \(\nabla u=\rho\nabla r\) is bounded there.

The inequalities \(H(r)\geq cI\) and \(\det H(\psi)=g_\theta\) give
\[
  \mathcal Lr\geq c\,\operatorname{tr}\big(H(\psi)^{-1}\big),\qquad
  \operatorname{tr}\big(H(\psi)^{-1}\big)\geq m g_\theta^{-1/m}\gtrsim 1 .
\]
The trace estimate is the arithmetic--geometric mean inequality
for the positive eigenvalues of the inverse Hessian. The Hessian
commutes with differentiation in any constant real unit direction
\(e\). Differentiating the equation once gives
\(\mathcal L\psi_e=(\log g_\theta)_e\). Differentiating a second time gives
\[
  (\log g_\theta)_{ee}
  = \mathcal L\psi_{ee}
    +D^2(\log\det)\big|_{H(\psi)}[H(\psi_e),H(\psi_e)]
  \leq \mathcal L\psi_{ee},
\]
by concavity of \(\log\det\). The derivatives of \(\log g_\theta\) are
uniformly bounded. Hence for a large fixed \(C\),
\(\mathcal L(\pm\psi_e+Cr)\geq0\) and \(\mathcal L(\psi_{ee}+Cr)\geq0\).
The maximum principle bounds these functions by their boundary
maxima. The first application gives a gradient bound throughout
the domain. The second reduces the upper Hessian bound to the
boundary.

It remains to bound the real Hessian on the boundary, uniformly in
\(\theta\). The next two steps treat the tangential and mixed entries
first, and then recover the last normal entry from the determinant.

\step{Tangential and mixed second derivatives.}
Fix a boundary point at the origin of a flat chart, and write
\(\mathbf n\) for its real outward unit normal and \(\nu\) for its
quaternionic column. If \(e,w\) are real tangent vectors at the
origin, differentiating \(u=\rho r\) gives \(D^2u(e,w)=\rho D^2r(e,w)\).
The derivatives of \(\rho\) are multiplied by \(r\) or by a tangential
derivative of \(r\), both of which vanish there.

To control a mixed entry \(D^2u(\mathbf n,e)\), we adapt the tangential
differentiation argument of Caffarelli, Nirenberg and Spruck
\cite[Section~2, (2.11)--(2.14)]{caffarelli-nirenberg-spruck}.
We want a rigid motion whose velocity at the origin is \(e\) and
whose normal component vanishes to first order along the boundary.
Let \(v\) be the real tangential projection of \(D^2r\,e/|\nabla r|\) at
the origin, viewed as a quaternionic column, and set
\[
  A=v\nu^*-\nu v^*-\nu(\nu^*v)\nu^* .
\]
The scalar \(\nu^*v\) is purely imaginary because \(v\) is real tangent.
Taking the adjoint and multiplying by \(\nu\) give \(A^*=-A\) and \(A\nu=v\).
The real vector field \(\mathcal V=(e+Ax)\cdot\nabla\) has
\(\mathcal Vr(0)=0\), and for any real tangent vector \(w\),
\[
  \partial_w(\mathcal Vr)(0)
   = D^2r(w,e)+|\nabla r|\,\mathbf n\cdot Aw
   = D^2r(w,e)-|\nabla r|\,w\cdot v=0 .
\]
The last equality follows from the definition of \(v\). Taylor's
theorem on the boundary gives \(\mathcal Vr=O(|x|^2)\) there, and
hence \(\mathcal Vu=\rho\mathcal Vr=O(|x|^2)\) there.

The local flow of \(\mathcal V\) consists of translations and quaternionic
unitary rotations. By covariance in Lemma~\ref{lem:hessian}, pulling
\(\psi\) back along this flow conjugates its Hessian by a unitary matrix
at the translated point. Its determinant is thus the pullback of
\(g_\theta\). Differentiating at flow time zero gives
\(\mathcal L(\mathcal V\psi)=\mathcal V\log g_\theta\). Subtracting the known
\(\mathcal L(\mathcal V B_0)\) gives
\(|\mathcal L(\mathcal V u)|\lesssim\operatorname{tr}(H(\psi)^{-1})\).

On the intersection of the domain with a small chart ball, take
the barrier \(C_1|x|^2-C_2r\). Choose \(C_1\) to dominate
\(|\mathcal Vu|\) on the true boundary by the quadratic bound and on
the sphere by the gradient bound. In four slots \(H(|x|^2)=8I\).
For \(C_2\) large, the lower bound on \(\mathcal Lr\) makes
\(\mathcal L(C_1|x|^2-C_2r\pm\mathcal Vu)\leq0\). Both functions are
nonnegative on the two boundary pieces, so they are nonnegative
inside by the minimum principle. Both vanish at the origin.
Taking inward derivatives there bounds
\(|\partial_{\mathbf n}\mathcal Vu(0)|\). Finally,
\[
  \partial_{\mathbf n}\mathcal Vu(0)
  =D^2u(\mathbf n,e)+\nabla u(0)\cdot A\nu
  =D^2u(\mathbf n,e),
\]
since \(\nabla u\) is normal and \(A\nu=v\) is real tangent.
All chart radii and constants can be chosen uniformly over the
compact boundary.

\step{The last normal derivative.}
The real tangent space is orthogonal to \(\mathbf n\) as real vectors.
To use the determinant we compress to the quaternionic orthogonal
complement \(\nu^\perp=\{w:\nu^*w=0\}\). On \(r=0\), the product rule gives
\[
  H(u)=\rho H(r)+G(\rho)G(r)^*+G(r)G(\rho)^*,
  \qquad G(r)=|\nabla r|\nu .
\]
The last two terms vanish on compression to \(\nu^\perp\), as in
\eqref{eq:boundary-compression}. Thus the compression of \(H(\psi)\)
is bounded below by \(b'H(r)|_{\nu^\perp}\gtrsim I\). The one real
Hessian entry of \(\psi\) we still need is \(\psi_{\mathbf n\mathbf n}\).
Its contribution to \(H(\psi)\) is \(\psi_{\mathbf n\mathbf n}\nu\nu^*\)
by the rank-one identity. In a quaternionic orthonormal basis
beginning with \(\nu\), the Hessian therefore takes the form
\[
  H(\psi)=
    \begin{pmatrix}\psi_{\mathbf n\mathbf n}+b_1&c^*\\c&D\end{pmatrix}.
\]
Here \(b_1,c,D\) involve only entries of the real Hessian that are
already bounded, and \(D\gtrsim I\). When \(m>1\), the Schur determinant
formula in Lemma~\ref{lem:cofactor} gives
\[
  \psi_{\mathbf n\mathbf n}
    =\frac{g_\theta}{\det D}+c^*D^{-1}c-b_1 .
\]
Every term on the right is bounded. When \(m=1\), the equation reads
\(\psi_{\mathbf n\mathbf n}+b_1=g_\theta\). The earlier maximum-principle
argument now bounds every unit directional second derivative from
above throughout the domain. Since \(\Delta\psi=\operatorname{tr}H(\psi)>0\),
each real Hessian eigenvalue is bounded below by minus the sum of
the upper bounds for the other eigenvalues. This yields the full
\(C^2\) bound along the path.

\step{Closedness and smooth dependence.}
The \(C^2\) estimate bounds the eigenvalues of \(H(\psi)\) from above.
They are positive and their product is \(g_\theta\), uniformly bounded
away from zero. Choose \(0<\lambda_-<\lambda_+\) so that, along the entire path, all these
eigenvalues lie strictly inside \([\lambda_-,\lambda_+]\), with a fixed margin at either end.

To apply the Evans--Krylov theory we extend the scalar operator, which so
far is defined only where the encoded Hessian is positive, to all real
symmetric matrices. In a flat chart
of the four-slot lift let \(S\) be any real symmetric \(4m\times4m\) matrix.
The encoding in Lemma~\ref{lem:hessian} is the linear map
\[
  \mathcal H:\operatorname{Sym}(4m,\mathbb R)\longrightarrow
  \operatorname{Herm}(m,\mathbb H),\qquad
  \mathcal H(S)_{ab}=\sum_{\alpha,\beta=0}^3
  e_\alpha\overline{e_\beta}\,S_{a\alpha,b\beta},
\]
so that \(H(\psi)=\mathcal H(D^2\psi)\). The map \(\mathcal H\) is far from
injective: when \(m=1\) it sends a real \(4\times4\) matrix to its trace, so it
forgets every trace-free direction. We never need to invert it. On the cone
\(\mathcal H(S)>0\) we can simply compose, setting
\(F(S)=\log\det\mathcal H(S)\), where the determinant is the Moore
determinant of the encoded matrix.

Write \(\mathcal E=\{Q\in\operatorname{Herm}(m,\mathbb H):\lambda_- I\leq Q\leq\lambda_+ I\}\).
At each \(Q\in \mathcal E\), take the tangent plane to \(\log\det\) and evaluate it
at \(\mathcal H(S)\). This gives a real affine function of \(S\), even when
\(\mathcal H(S)\) is not positive. Define
\[
  \ell_Q(S)=\log\det Q+Q^{-1}:\big(\mathcal H(S)-Q\big),
  \qquad \widetilde F(S)=\inf_{Q\in \mathcal E}\ell_Q(S).
\]
The infimum is finite since \(\mathcal E\) is compact. Concavity puts each tangent
plane \emph{above} the graph of \(\log\det\) on its positive cone. When
\(\mathcal H(S)\in \mathcal E\), the choice \(Q=\mathcal H(S)\) gives equality, so
\(\widetilde F=F\) there. On all real symmetric matrices \(\widetilde F\)
is concave, being an infimum of affine functions; compare the construction
in~\cite[Lemma~4.1]{tosatti-et-al}.

The kernel of \(\mathcal H\) does not spoil ellipticity, because real
ellipticity only tests \emph{positive} increments. For a rank-one increment Lemma~\ref{lem:hessian} gives
\[
  \ell_Q(S+\xi\xi^{\mathsf t})-\ell_Q(S)
   =Q^{-1}:\widehat\xi\widehat\xi^*
   =\widehat\xi^* Q^{-1}\widehat\xi .
\]
Since \(|\widehat\xi|=|\xi|\), this lies between \(|\xi|^2/\lambda_+\) and
\(|\xi|^2/\lambda_-\), independently of \(Q\). Decompose any real \(P\geq0\) into
rank-one increments and add. The resulting bounds hold for every
\(\ell_Q\) with the same constants, so taking the infimum preserves them:
\[
  \frac{1}{\lambda_+}\operatorname{tr}P
  \leq \widetilde F(S+P)-\widetilde F(S)
  \leq \frac{1}{\lambda_-}\operatorname{tr}P \qquad(P\geq0).
\]
This is uniform ellipticity in the full real matrix space. Moreover,
\(\widetilde F(D^2\psi)=\log g_\theta\) along the path, so the extension
has not changed the equation at any solution.

The right-hand sides \(\log g_\theta\) have uniform smooth bounds and the
Dirichlet values are the fixed smooth function \(B_0\). We may now use the
interior Evans--Krylov estimate for a concave operator
\cite[Theorem~3.4]{tosatti-et-al} and the boundary \(C^{2,\alpha}\) estimate
\cite[Theorem~1.2]{silvestre-sirakov}. These versions do not require
the envelope to be twice differentiable away from the solution Hessians.
Combining the boundary expansions with the interior estimates on balls
of radius comparable to the distance to the boundary, then using a finite
chart cover, gives a uniform \(C^{2,\alpha}\) bound for some \(\alpha>0\).
For a sequence of solvable parameters tending to a limit, a
subsequence converges in \(C^{2,\beta}\) for \(0<\beta<\alpha\).
The equation, boundary values and strict positivity pass to the limit.

For the bootstrap we return to the original \(F\), which is smooth on the
open cone containing the solution Hessians. Flatten the boundary and apply Schauder
estimates to the tangential derivatives with their prescribed boundary
values. These estimates gain the additional derivatives involving
a tangential direction. Uniform ellipticity lets us solve the equation
locally for the pure normal second derivative as a smooth function
of the other entries and lower derivatives. Iterating gives smoothness
up to the boundary. Thus the solvable set is closed and reaches
\(\theta=1\).

Translations in the added torus preserve the data and the equation.
Uniqueness makes the lifted solution invariant under these translations,
so it descends to the three-slot domain. For a smooth family of moving
domains, use its local identification with a fixed domain to pull the
lifted equation back. The coefficients and boundary data depend
smoothly on the parameters. The linearization is conjugate to the
Dirichlet isomorphism used for openness. The implicit function
theorem in successively higher Hölder spaces gives smooth dependence
up to the boundary. The strict determinant inequalities used to
choose \(b,b'\) persist for nearby parameters, so the same constants
work locally.
\end{proof}

\end{document}